\documentclass[12pt, reqno]{amsart}
\usepackage{eurosym}
\usepackage{amsmath, amsthm, amscd, amsfonts, amssymb, graphicx, color}
\usepackage[bookmarksnumbered, colorlinks, plainpages]{hyperref}

\hypersetup{colorlinks=true,linkcolor=red, anchorcolor=green, citecolor=cyan, urlcolor=red, filecolor=magenta, pdftoolbar=true}
\newtheorem{theorem}{Theorem}[section]

\newtheorem{proposition}[theorem]{Proposition}
\newtheorem{corollary}[theorem]{Corollary}
\theoremstyle{definition}
\newtheorem{definition}[theorem]{Definition}
\newtheorem{example}[theorem]{Example}

\theoremstyle{remark}

\numberwithin{equation}{section}

\begin{document}
\title[Composition ideals of Lip-Linear operators]{Composition ideals of
Lip-Linear operators and a Hilbert space characterization}
\author{ Athmane Ferradi and Khalil Saadi}
\date{}
\dedicatory{Laboratory of Functional Analysis and Geometry of Spaces,
Faculty of Mathematics and Computer Science, University of Mohamed
Boudiaf-M'sila, Po Box 166, Ichebilia, 28000, M'sila, Algeria.\\
Second address for the first author: Ecole Normale Sup\'{e}rieure de
Bousaada, Algeria.\\
athmane.ferradi@univ-msila.dz\\
khalil.saadi@univ-msila.dz}

\begin{abstract}
In this paper, we investigate classes of Lip-linear operators constructed
using the composition ideal method. We focus on two fundamental linear
operator ideals, $p$-summing and strongly $p$-summing operators, and extend
them to define the corresponding classes of Lip-linear operators. Several
key results are established, including a characterization theorem for
Hilbert spaces originally due to Kwapie\'{n}. Specifically, we show that a
Banach space $F$ is isomorphic to a Hilbert space if and only if every
factorable strongly $p$-summing Lip-linear operator with values in $F$ is
Cohen strongly $p$-summing.
\end{abstract}

\maketitle

\setcounter{page}{1}


\let\thefootnote\relax\footnote{
\par
{}\textit{Key words and phrases. }Lip-Linear operators; Cohen strongly $p$%
-summing Lip-Linear operators; strongly $p$-summing Lip-Linear operators;
linear ideals; non linear ideals, Pietsch factorization; Kwapie\'{n} theorem.%
}

\section{Introduction and preliminaries}

The theory of operator ideals has been a central theme in the study of
Banach spaces, with deep implications for linear and nonlinear analysis.
Motivated by the success of this theory in the linear setting, many efforts
have been made to extend these ideas to more general classes of mappings,
including multilinear operators, homogeneous polynomials, and recently,
Lipschitz mappings. Composition ideals are among the key techniques used to
construct new classes of nonlinear operator ideals. In the multilinear
setting, G. Botelho et al. \cite{PelCompo} studied and developed the theory
of composition ideals for multilinear mappings and homogeneous polynomials.
Given an operator ideal $\mathcal{I}$, they defined the associated
multi-ideal $\mathcal{I}\circ \mathcal{L}$ as the class of multilinear
mappings that factor through an operator in $\mathcal{I}$. This construction
preserves many of the structural properties of the original ideal and
provides a bridge between linear and multilinear operator theory. In \cite%
{FS2025}, the class of Lip-linear operators has been introduced and
staudied. These operators are Lipschitz with respect to the first component
and linear in the second. Our goal is to study operator ideals in the
context of Lip-linear operators. We define and investigate, in particulier,
composition ideals of Lip-linear operators. Let $X$ be a pointed metric
space, and let $E,F$ be Banach spaces. A mapping $T:X\times E\rightarrow F$
is called Lip-linear if it is Lipschitz in the first variable and linear in
the second. The space of all such operators $T$ satisfying $T\left(
0,e\right) =0$ for every $e\in E$ is denoted by $LipL_{0}(X\times E;F).$
This space becomes a Banach space when equipped with the norm%
\begin{equation*}
LipL\left( T\right) =\sup_{x\neq y,\left\Vert e\right\Vert =1}\frac{%
\left\Vert T\left( x,e\right) -T\left( y,e\right) \right\Vert }{d\left(
x,y\right) }.
\end{equation*}%
These operators exhibit strong structural properties, particularly through
their connection with the Lipschitz tensor product $X\boxtimes E$. Every
Lip-linear operator $T:X\times E\rightarrow F$ admits a unique linearization 
$\widehat{T}:X\widehat{\boxtimes }_{\pi }E\rightarrow F$ defined by%
\begin{equation*}
\widehat{T}(\sum\limits_{i=1}^{n}\delta _{\left( x_{i},y_{i}\right)
}\boxtimes e_{i})=\sum\limits_{i=1}^{n}T\left( x_{i},e_{i}\right) -T\left(
y_{i},e_{i}\right) ,
\end{equation*}%
and conversely. This correspondence establishes an isometric identification
between the class $LipL_{0}(X\times E;F)$ and the space of bounded linear
operators $\mathcal{L}(X\widehat{\boxtimes }_{\pi }E;F)$. Composing linear
operators with Lip-linear operators yields new classes of Lip-linear
operators. The associated linearized operator characterizes this relation
and preserves the structure through the composition. This approach, known as
the composition method, provides a systematic way to define new classes of
Lip-linear operators. As in the multilinear case, we use this method to
construct and analyze specific examples of Lip-linear operator ideals.

This paper is structured as follows: Section 1 introduces the essential
definitions and foundational concepts used throughout the article, including
Lipschitz mappings, linear and bilinear operators, as well as key
identification results. It also recalls several summability notions for
Lipschitz, linear and bilinear operators. In Section 2, we define Lip-linear
operator ideals as a natural generalization of multilinear ideals. We focus
in particular on those constructed using the composition method. Several
structural results are established in this setting. Section 3 is devoted to
the introduction of Cohen strongly $p$-summing Lip-linear operators. We show
that this class forms a Lip-linear ideal, interpreted as the composition of
the ideal of strongly $p$-summing linear operators $\mathcal{D}_{p}$ with
the space of Lip-linear operators. In Section 4, motivated by the
characterization of Hilbert spaces due to Kwapie\'{n} \cite{Kwa}, we
introduce the notion of strongly $p$-summing Lip-linear operators as a
natural extension of Dimant's concept in the multilinear case \cite{Dim}.
However, since these operators do not, in general, correspond to their
linearizations under the usual concept of $p$-summing perators, we propose a
refined class: the factorable strongly $p$-summing Lip-linear operators.
This refined notion allows us to recover the connection with the associated
linearized operators. One of the main results of this paper is a
characterization of Hilbert spaces using Lip-linear operators. We prove that
a Banach space $F$ is isomorphic to a Hilbert space if and only if every
factorable strongly $p$-summing Lip-linear operator with values in $F$ is
Cohen strongly $p$-summing.

The letters $E,F$ and $G$ stand for Banach spaces over $\mathbb{R}$. The
closed unit ball of $E$ is denoted by $B_{E}$ and its topological dual by $%
E^{\ast }.$ The Banach space of all continuous linear operators from $E$ to $%
F$ is denoted by $\mathcal{L}(E;F)$, and $\mathcal{B}(E\times F;G)$ denotes
the Banach space of all continuous bilinear operators, equipped with the
standard sup norm. Let $X$ and $Y$ be pointed metric spaces with base points
denoted by $0$. The Banach space of Lipschitz functions $T:X\rightarrow E$
that satisfy $T\left( 0\right) =0$ is denoted $Lip_{0}(X;E)$ and equipped
with the norm%
\begin{equation*}
Lip(T)=\sup_{x\neq y}\frac{\left\Vert T\left( x\right) -T\left( y\right)
\right\Vert }{d\left( x,y\right) }.
\end{equation*}%
The closed unit ball $B_{X^{\#}}$ of $X^{\#}:=Lip_{0}(X;\mathbb{K})$ forms a
compact Hausdorff space under the topology of pointwise convergence on $X$ .
Define the mapping $\delta _{X}:X\longrightarrow \left( X^{\#}\right) ^{\ast
}$ by $\delta _{x}(f)=f\left( x\right) $ for all $f\in X^{\#}.$ For $x,y\in X
$ and $f\in X^{\#}$, we have%
\begin{equation*}
\delta _{\left( x,y\right) }\left( f\right) =\delta _{x}(f)-\delta
_{y}(f)=f\left( x\right) -f\left( y\right) .
\end{equation*}%
We have $\delta _{\left( x,y\right) }=\delta _{\left( x,0\right) }-\delta
_{\left( y,0\right) }.$ The closed subspace spanned by $\left\{ \delta
_{\left( x,y\right) }:x,y\in X\right\} $ is a Banach space, denoted by $%
\mathcal{F}(X),$ with the following norm%
\begin{equation*}
\Vert m\Vert _{\mathcal{F}(X)}=\inf \left\{ \sum_{i=1}^{n}|\lambda
_{i}|d(x_{i},y_{i}):\ m=\sum_{i=1}^{n}\lambda _{i}\delta _{\left(
x_{i}y_{i}\right) }\right\} 
\end{equation*}%
where the infimum is taken over all representations of $m$. It is well known
that $\mathcal{F}(X)^{\ast }=X^{\#}.$ For any $T\in Lip_{0}(X;E)$, there is
a unique linearization $T_{L}\in \mathcal{L}(\mathcal{F}(X);E)$ such that%
\begin{equation*}
T=T_{L}\circ \delta _{X},
\end{equation*}%
with $\Vert T_{L}\Vert =Lip(T)$. This yields the identification%
\begin{equation}
Lip_{0}(X;E)=\mathcal{L}(\mathcal{F}(X);E).  \label{1.0}
\end{equation}%
If $T:E\rightarrow F$ is a linear operator between Banach spaces, then the
corresponding linear operator $T_{L}$ is given by%
\begin{equation*}
T_{L}=T\circ \beta _{E},
\end{equation*}%
where $\beta _{E}:\mathcal{F}\left( E\right) \rightarrow E$ is linear
quotient map verifying $\beta _{E}\circ \delta _{E}=id_{E}$ and $\left\Vert
\beta _{E}\right\Vert \leq 1$ (see \cite[Lemma 2.4]{gk} for more details
about $\beta _{E}$). We denote by $X\boxtimes E$ the Lipschitz tensor
product of $X$ and $E,$ as introduced and studied by Cabrera-Padilla et al.
in \cite{ccjv}. This space is spanned by the functionals $\delta _{\left(
x,y\right) }\boxtimes e$ on $Lip_{0}(X;E^{\ast }),$ defined by $\delta
_{\left( x,y\right) }\boxtimes e(f)=\langle f(x)-f(y),e\rangle $ for $x,y\in
X$ and $e\in E.$ The Lipschitz projective norm on $X\boxtimes E$ is given by 
\begin{equation*}
\pi (u)=\inf \left\{ \sum_{i=1}^{n}d(x_{i},y_{i})\Vert e_{i}\Vert ,\
u=\sum_{i=1}^{n}\delta _{\left( x_{i},y_{i}\right) }\boxtimes e_{i}\right\} ,
\end{equation*}%
where the infimum is taken over all representations of $u$ and $X\widehat{%
\boxtimes }_{\pi }E$ denotes its completion. The Banach space $%
Lip_{0}(X;E^{\ast })$ is isometrically isomorphic to $(X\widehat{\boxtimes }%
_{\pi }E)^{\ast }$. According to \cite[Proposition 6.7]{ccjv}, $X\widehat{%
\boxtimes }_{\pi }E$ is isometrically isomorphic to $\mathcal{F}\left(
X\right) \widehat{\otimes }_{\pi }E$ via the mapping 
\begin{equation*}
I\left( u\right) =I(\sum\limits_{i=1}^{n}\delta _{\left( x_{i};y_{i}\right)
}\boxtimes e_{i})=\sum\limits_{i=1}^{n}\delta _{\left( x_{i};y_{i}\right)
}\otimes e_{i}.
\end{equation*}%
This leads to the following isometric identifications%
\begin{equation*}
\mathcal{B}\left( \mathcal{F}\left( X\right) \times E;F\right) =\mathcal{L}%
\left( \mathcal{F}\left( X\right) \widehat{\otimes }_{\pi }E;F\right) =%
\mathcal{L}\left( X\widehat{\boxtimes }_{\pi }E;F\right) .
\end{equation*}%
A mapping $T:X\times E\rightarrow F$ is called a Lip-linear operator if
there is a constant $C>0$ such that for all $x,y\in X$ and $e\in E$, 
\begin{equation}
\left\Vert T(x,e)-T(y,e)\right\Vert \leq Cd(x,y)\left\Vert e\right\Vert ,
\label{1.1}
\end{equation}%
and $T(x,\cdot )$ is linear for each fixed $x\in X$. By definition, $T(x,0)=0
$ for every $x\in X$. We denote by $LipL_{0}(X\times E;F)$ the space of all
Lip-linear operators from $X\times E$ to $F$ satisfying $T(0,e)=0$ for every 
$e\in E.$ This space becomes a Banach space when equipped with the norm 
\begin{equation*}
LipL\left( T\right) :=\inf \left\{ C:\text{ }C\text{ satisfies }\left( \ref%
{1.1}\right) \right\} (=\sup_{x\neq y\text{ }\left\Vert e\right\Vert =1}%
\frac{\left\Vert T\left( x,e\right) -T\left( y,e\right) \right\Vert }{%
d\left( x,y\right) }).
\end{equation*}%
For a Lip-linear operator $T:X\times E\rightarrow F,$ the following diagram
is commutative%
\begin{equation}
\begin{array}{ccc}
X\times E & \overset{T}{\longrightarrow } & F \\ 
\sigma ^{LipL}\downarrow  & \nearrow \widehat{T} &  \\ 
X\boxtimes E &  & 
\end{array}
\label{1.1.2}
\end{equation}%
where $\sigma ^{LipL}=\delta _{X}\boxtimes id_{E},$ and $LipL\left( \sigma
^{LipL}\right) =1.$ When $F=\mathbb{R}$ we denote the space simply by $%
LipL_{0}(X\times E).$ Let $T:X\times E\longrightarrow F$ be a Lip-linear
operator. By \cite{FS2025}, there exists a unique bilinear operator $T_{bil}:%
\mathcal{F}\left( X\right) \times E\rightarrow F$ such that 
\begin{equation}
T=T_{bil}\left( \delta _{X},id_{E}\right) ,  \label{1.2}
\end{equation}
and the identification 
\begin{equation}
LipL_{0}(X\times E;F)=\mathcal{B}(\mathcal{F}\left( X\right) \times E;F),
\label{1.4}
\end{equation}%
holds isometrically. Moreover, we have the following commutative diagram 
\begin{equation}
\begin{array}{ccc}
\mathcal{F}\left( X\right) \times E & \overset{T_{bil}}{\longrightarrow } & F
\\ 
\sigma _{2}\downarrow  & \nearrow \widetilde{T_{bil}} & \uparrow \widehat{T}
\\ 
\mathcal{F}\left( X\right) \widehat{\otimes }_{\pi }E & \overset{I^{-1}}{%
\longrightarrow } & X\widehat{\boxtimes }_{\pi }E%
\end{array}
\label{1.5}
\end{equation}%
where $\sigma _{2}$ is the canonical bilinear map defined by $\left(
m,e\right) \mapsto m\otimes e$ and $\widetilde{T_{bil}}$ is the
linearization of $T_{bil}.$

Let us recall some fundamental concepts: Let $X$ be a pointed metric space
and $E,F$ be Banach spaces.

- The linear operator $u$ is called $p$-\textit{summing, }$u\in \left( \Pi
_{p}(E;F),\pi _{p}\left( \cdot \right) \right) ,$ if there exists a constant 
$C>0$ such that for all $\left( x_{i}\right) _{i=1}^{n}\subset E$ we have:%
\begin{equation*}
(\sum\limits_{i=1}^{n}\left\Vert u\left( x_{i}\right) \right\Vert ^{p})^{%
\frac{1}{p}}\leq C\left\Vert \left( x_{i}\right) \right\Vert _{\ell
_{p}^{n,w}\left( E\right) }.
\end{equation*}

- The linear operator $u$ is called \textit{strongly }$p$-\textit{summing, }$%
u\in \left( \mathcal{D}_{p}(E;F),d_{p}\left( \cdot \right) \right) ,$ if
there exists a constant $C>0$ such that for all $\left( x_{i}\right)
_{i=1}^{n}\subset E$ and $\left( z_{i}^{\ast }\right) _{i=1}^{n}\subset
F^{\ast }$ we have:%
\begin{equation*}
\sum\limits_{i=1}^{n}\left\vert \left\langle u\left( x_{i}\right)
,z_{i}^{\ast }\right\rangle \right\vert \leq C\left\Vert \left( x_{i}\right)
\right\Vert _{\ell _{p}^{n}\left( E\right) }\left\Vert \left( z_{i}^{\ast
}\right) \right\Vert _{\ell _{p^{\ast }}^{n,w}\left( F^{\ast }\right) }.
\end{equation*}

- The Lipschitz operator $R$ is called \textit{Lipschitz strongly }$p$-%
\textit{summing, }$u\in (\mathcal{D}_{p}^{L}(X;F),d_{p}^{L}\left( \cdot
\right) ),$ if there exists a constant $C>0$ such that for all $\left(
x_{i}\right) _{i=1}^{n},\left( y_{i}\right) _{i=1}^{n}\subset X$ and $\left(
z_{i}^{\ast }\right) _{i=1}^{n}\subset F^{\ast }$ we have:%
\begin{equation*}
\sum\limits_{i=1}^{n}\left\vert \left\langle R\left( x_{i}\right) -R\left(
y_{i}\right) ,z_{i}^{\ast }\right\rangle \right\vert \leq
C(\sum\limits_{i=1}^{n}d(x_{i},y_{i})^{p})^{\frac{1}{p}}\left\Vert \left(
z_{i}^{\ast }\right) \right\Vert _{\ell _{p^{\ast }}^{n,w}\left( F^{\ast
}\right) }.
\end{equation*}

Noted that, the nomrs of above classes are defined as the infimum of all
constants $C$ satisfying the respective inequalities.

\begin{definition}
A bilinear operator $B:E\times F\longrightarrow G$\ is Cohen strongly $p$%
\textit{-summing}\ if there exists a constant $C>0$\ such that, for any $%
\left( x_{i}\right) _{i=1}^{n}\subset E,\left( y_{i}\right)
_{i=1}^{n}\subset F$ and $\left( z_{i}^{\ast }\right) _{i=1}^{n}\subset
G^{\ast }$, we have 
\begin{equation}
\sum\limits_{i=1}^{n}\left\vert \left\langle B(x_{i},y_{i}),z_{i}^{\ast
}\right\rangle \right\vert \leq C(\sum\limits_{i=1}^{n}\left( \left\Vert
x_{i}\right\Vert \left\Vert y_{i}\right\Vert \right) ^{p})^{\frac{1}{p}%
}\left\Vert \left( z_{i}^{\ast }\right) \right\Vert _{\ell _{p^{\ast
}}^{n,w}\left( F^{\ast }\right) }.  \label{1.3}
\end{equation}
\end{definition}

We denote by $\mathcal{D}_{p}^{2}(E\times F;G)$ the space of all Cohen
strongly $p$-smming bilinear operators from $E\times F$ to $G$. For $B\in 
\mathcal{D}_{p}^{2}(E\times F;G)$, we define 
\begin{equation*}
d_{p}^{2}(B)=\inf \left\{ C:\text{ }C\text{ satisfies }\left( \ref{1.3}%
\right) \right\} .
\end{equation*}

\section{Composition method of Lip-linear Operators}

Several studies have investigated linear and nonlinear operator ideals (see 
\cite{atlb}, \cite{gbet}, \cite{PelCompo}, \cite{dk}, \cite{sgh}). In this
section, we extend this theory to define Lip-linear operator ideals. Using
the composition method, we construct Lip-linear ideals from existing linear
operator ideals. This approach broadens the theory and enables the transfer
of many properties from the linear to the Lip-linear case.

\begin{definition}
\label{Defi1}A Lip-Linear ideal $\mathcal{I}^{LipL}$, is a subclass of $%
LipL_{0}$ such that for every pointed metric space $X$ and every Banach
spaces $E$ and $F$ the components 
\begin{equation*}
\mathcal{I}^{LipL}(X\times E;F):=LipL_{0}(X\times E;F)\cap \mathcal{I}%
^{LipL},
\end{equation*}%
satisfy:

\begin{enumerate}
\item[(i)] $\mathcal{I}^{LipL}(X\times E;F)$ is a vector subspace of $%
LipL_{0}(X\times E;F)$.

\item[(ii)] The map $f\cdot e^{\ast }\cdot z$ belongs to $\mathcal{I}%
^{LipL}(X\times E;F)$, for every $f\in X^{\#},$ $e^{\ast }\in E^{\ast }$ and 
$z\in F$.

\item[(iii)] Ideal property: If $f\in Lip_{0}(X;X_{0})$, $v\in \mathcal{L}%
(E;E_{0})$, $T\in \mathcal{I}^{LipL}(X_{0}\times E_{0};F_{0})$ and $u\in 
\mathcal{L}(F_{0};F)$ then $u\circ T\circ (f,v)\in \mathcal{I}%
^{LipL}(X\times E;F)$.

A Lip-Linear ideal $\mathcal{I}^{LipL}$ is a Banach Lip-Linear ideal if
there is $\Vert \cdot \Vert _{\mathcal{I}^{LipL}}:\mathcal{I}%
^{LipL}\longrightarrow \left[ 0,+\infty \right[ $ that satisfies:

\item[(i')] For every pointed metric space $X$ and every Banach spaces $E$
and $F$, the pair $(\mathcal{I}^{LipL}(X\times E;F),\Vert \cdot \Vert _{%
\mathcal{I}^{LipL}})$ is a Banach space and $LipL\left( T\right) \leq \Vert
T\Vert _{\mathcal{I}^{LipL}}$ for all $T\in \mathcal{I}^{LipL}\left( X\times
E;F\right) $.

\item[(ii')] $\Vert Id_{\mathbb{R}^{2}}:\mathbb{R}\times \mathbb{R}%
\longrightarrow \mathbb{R}:Id_{\mathbb{R}^{2}}(\alpha ,\beta )=\alpha \beta
\Vert _{\mathcal{I}^{LipL}}=1$.

\item[(iii')] If $f\in Lip_{0}(X;X_{0})$, $v\in \mathcal{L}(E;E_{0})$, $T\in 
\mathcal{I}^{LipL}(X_{0}\times E_{0};F_{0})$ and $u\in \mathcal{L}(F_{0};F)$%
, the inequality 
\begin{equation*}
\left\Vert u\circ T\circ (f,v)\right\Vert _{\mathcal{I}^{LipL}}\leq \Vert
u\Vert \left\Vert T\right\Vert _{\mathcal{I}^{LipL}}Lip(f)\Vert v\Vert 
\end{equation*}%
holds.
\end{enumerate}
\end{definition}

Let $T$ be a Lip-linear operator. We define its dual (transpose) $%
T^{t}:F^{\ast }\rightarrow LipL_{0}\left( X\times E\right) $ by%
\begin{equation*}
T^{t}\left( y^{\ast }\right) \left( x,e\right) =\left\langle y^{\ast
},T\left( x,e\right) \right\rangle ,
\end{equation*}%
for all $y^{\ast }\in F^{\ast },x\in X$ and $e\in E.$ We have%
\begin{eqnarray*}
LipL\left( T\right)  &=&\sup_{x\neq y,\left\Vert e^{\ast }\right\Vert =1}%
\frac{\left\Vert T\left( x,e\right) -T\left( y,e\right) \right\Vert }{%
d\left( x,y\right) } \\
&=&\sup_{x\neq y,\left\Vert e^{\ast }\right\Vert =1}\sup_{\left\Vert z^{\ast
}\right\Vert =1}\frac{\left\vert z^{\ast }\left( T\left( x,e\right) \right)
)-z^{\ast }\left( T\left( y,e\right) \right) \right\vert }{d\left(
x,y\right) } \\
&=&\sup_{x\neq y,\left\Vert e^{\ast }\right\Vert =1}\sup_{\left\Vert z^{\ast
}\right\Vert =1}\frac{\left\vert T^{t}\left( z^{\ast }\right) \left(
x,e\right) -T^{t}\left( z^{\ast }\right) \left( y,e\right) \right\vert }{%
d\left( x,y\right) } \\
&=&\sup_{\left\Vert z^{\ast }\right\Vert =1}LipL\left( T^{t}\left( z^{\ast
}\right) \right) =\left\Vert T^{t}\right\Vert .
\end{eqnarray*}%
Consider a simple construction of a Lip-linear ideal based on a given
operator ideal $\mathcal{I}$. We define a corresponding class of Lip-linear
operators by%
\begin{equation*}
\mathcal{I}_{dual}^{LipL}(X\times E;F)=\left\{ T\in LipL_{0}(X\times
E;F):T^{t}\in \mathcal{I}(F^{\ast };LipL_{0}(X\times E)\right\} ,
\end{equation*}%
Furthermore, if $(\mathcal{I},\left\Vert .\right\Vert _{\mathcal{I}})$ is a
normed operator ideal, we equip this class with the norm%
\begin{equation}
\Vert T\Vert _{\mathcal{I}_{dual}^{LipL}}=\Vert T^{t}\Vert _{\mathcal{I}}.
\label{2.1}
\end{equation}

\begin{proposition}
\label{Propo1}Let $(\mathcal{I},\left\Vert \cdot \right\Vert _{\mathcal{I}})$
be a Banach operator ideal. Then, the class $(\mathcal{I}_{dual}^{LipL},%
\left\Vert \cdot \right\Vert _{\mathcal{I}_{dual}^{LipL}})$ is a Banach
Lip-Linear ideal.
\end{proposition}

\begin{proof}
We verify the ideal property. Let $f\in Lip_{0}(X;X_{0})$, $v\in \mathcal{L}%
(E;E_{0})$, $T\in \mathcal{I}_{dual}^{LipL}(X_{0}\times E_{0};F_{0})$ and $%
u\in \mathcal{L}(F_{0};F).$ We have%
\begin{equation*}
\left( u\circ T\circ (f,v)\right) ^{t}=\Psi _{\left( f,v\right) }\circ
T^{t}\circ u^{\ast },
\end{equation*}%
where the linear map $\Psi _{\left( f,v\right) }:LipL_{0}\left( X_{0}\times
E_{0}\right) \rightarrow LipL_{0}\left( X\times E\right) $ is defined by%
\begin{equation*}
\Psi _{\left( f,v\right) }\left( S\right) \left( x,e\right) =S\left( f\left(
x\right) ,v\left( e\right) \right) ,
\end{equation*}%
for all $S\in LipL_{0}\left( X_{0}\times E_{0}\right) .$ Clearly,%
\begin{equation*}
\Vert \Psi _{\left( f,v\right) }\Vert \leq Lip\left( f\right) \left\Vert
v\right\Vert .
\end{equation*}%
By the ideal property of $\mathcal{I}$, we conclude that 
\begin{equation*}
\left( u\circ T\circ (f,v)\right) ^{t}\in \mathcal{I}(F^{\ast
};LipL_{0}(X\times E).
\end{equation*}%
Moreover,%
\begin{eqnarray*}
\Vert u\circ T\circ (f,v)\Vert _{\mathcal{I}_{dual}^{LipL}} &=&\Vert \Psi
_{\left( f,v\right) }\circ T^{t}\circ u^{\ast }\Vert _{\mathcal{I}} \\
&\leq &\Vert \Psi _{\left( f,v\right) }\Vert \Vert T^{t}\Vert _{\mathcal{I}%
}\Vert u^{\ast }\Vert  \\
&\leq &Lip\left( f\right) \left\Vert v\right\Vert \Vert T\Vert _{\mathcal{I}%
_{dual}^{LipL}}\Vert u\Vert .
\end{eqnarray*}%
Finally, for every $T\in \mathcal{I}_{dual}^{LipL}(X\times E;F),$ we have%
\begin{equation*}
LipL\left( T\right) =\left\Vert T^{t}\right\Vert \leq \Vert T^{t}\Vert _{%
\mathcal{I}}=\Vert T\Vert _{\mathcal{I}_{dual}^{LipL}}.
\end{equation*}%
Hence, $\mathcal{I}_{dual}^{LipL}$ satisfies the ideal property and defines
a normed Lip-linear operator ideal.
\end{proof}

\begin{definition}
Let $\mathit{{\mathcal{I}}}$ be an operator ideal. An operator $T\in
LipL_{0}(X\times E;F)$ belongs to $\mathcal{I}\circ LipL_{0}$, we write $%
T\in \mathcal{I}\circ LipL_{0}(X\times E;F),$ if there are a Banach space $%
G, $ $S\in LipL_{0}(X\times E;G)$ and a linear operator $u\in \mathcal{I}%
(G;F)$ such that the following diagram commuts%
\begin{equation*}
\begin{array}{ccc}
X\times E & \overset{T}{\longrightarrow } & F \\ 
S\downarrow & \nearrow u &  \\ 
G &  & 
\end{array}%
\end{equation*}%
In the other hand, $T=u\circ S.$ If $(\mathcal{I},\left\Vert \cdot
\right\Vert _{\mathcal{I}})$ is a normed operator ideal, we write 
\begin{equation}
\Vert T\Vert _{\mathcal{I}\circ LipL_{0}}=\inf \{\Vert u\Vert _{\mathcal{I}%
}LipL(S):T=u\circ S\}.  \label{2.2}
\end{equation}
\end{definition}

\begin{theorem}
Let $\mathcal{I}$ be a Banach operator ideal. The pair $\left( \mathcal{I}%
\circ LipL_{0},\Vert \cdot \Vert _{\mathcal{I}\circ LipL_{0}}\right) $ forms
a Banach Lip-Linear ideal.
\end{theorem}

\begin{proof}
Following the same reasoning as in \cite[Proposition 3.7]{PelCompo}, we can
show that $\mathcal{I}\circ LipL_{0}$ constitutes a Banach Lip-Linear ideal.
Moreover, the norm in $\left( \ref{2.2}\right) $ is well-defined and
satisfies all the properties required by Definition $\ref{Defi1}$.
\end{proof}

\begin{proposition}
\label{CopoMethod}Let $(\mathcal{I},\Vert \cdot \Vert _{\mathcal{I}})$ be a
Banach operator ideal. An operator $T\in LipL_{0}(X\times E;F)$ belongs to $%
\mathcal{I}\circ LipL_{0}(X\times E;F),$ if and only if, its linearization $%
\widehat{T}$ belongs to $\mathcal{I}(X\widehat{\boxtimes }_{\pi }E;F).$
Furthermore, 
\begin{equation*}
\Vert T\Vert _{\mathcal{I}\circ LipL_{0}}=\Vert \widehat{T}\Vert _{\mathcal{I%
}},
\end{equation*}%
and we have the following isemotric identification%
\begin{equation*}
\mathcal{I}\circ LipL_{0}(X\times E;F)=\mathcal{I}(X\widehat{\boxtimes }%
_{\pi }E;F).
\end{equation*}
\end{proposition}

\begin{proof}
Suppose $\widehat{T}\in \mathcal{I}(X\widehat{\boxtimes }_{\pi }E;F)$.
Consider the factorization of $T$ given in $\left( \ref{1.1.2}\right) ,$ $T=%
\widehat{T}\circ \sigma ^{LipL}.$ Since $\sigma ^{LipL}$ is Lip-Linear
operator with $LipL(\sigma ^{LipL})=1,$ it follows that $T\in \mathcal{I}%
\circ LipL_{0}(X\times F;F).$ Moreover,%
\begin{equation*}
\left\Vert T\right\Vert _{\mathcal{I}\circ LipL_{0}}\leq \left\Vert \widehat{%
T}\right\Vert _{\mathcal{I}}.
\end{equation*}%
Conversely, asumme $T\in \mathcal{I}\circ LipL_{0}(X\times E;F)$. Then,
there exist a Banach space $G$, an operator $S\in LipL_{0}(X\times E;G)$ and
a linear operator $u\in \mathcal{I}(G;F)$ such that $T=u\circ S$. Since $%
\widehat{T}=u\circ \widehat{S},$ it follows from the ideal property that $%
\widehat{T}\in \mathcal{I}(X\widehat{\boxtimes }_{\pi }E;F)$. We have 
\begin{equation*}
\left\Vert \widehat{T}\right\Vert _{\mathcal{I}}\leq \left\Vert u\right\Vert
_{\mathcal{I}}\left\Vert \widehat{S}\right\Vert =\left\Vert u\right\Vert _{%
\mathcal{I}}LipL(S).
\end{equation*}%
Taking the infimum over all such factorizations yields%
\begin{equation*}
\left\Vert \widehat{T}\right\Vert _{\mathcal{I}}\leq \left\Vert T\right\Vert
_{\mathcal{I}\circ LipL_{0}}.
\end{equation*}%
To conclude, we verify surjectivity. Let $A\in \mathcal{I}(X\widehat{%
\boxtimes }_{\pi }E;F),$ we define $T$ by%
\begin{equation*}
T=A\circ \sigma ^{LipL}.
\end{equation*}%
We have $T\in \mathcal{I}\circ LipL_{0}(X\times F;F)$ and $\widehat{T}=A.$
\end{proof}

\begin{theorem}
If $\mathcal{I}$ is an operaor ideal, then 
\begin{equation*}
\mathcal{I}_{dual}^{LipL}=\mathcal{I}^{dual}\circ LipL_{0}.
\end{equation*}%
Moreover, if $(\mathcal{I},\left\Vert .\right\Vert _{\mathcal{I}})$ is a
Banach operator ideal, we have%
\begin{equation*}
\Vert \cdot \Vert _{\mathcal{I}_{dual}^{LipL}}=\Vert \cdot \Vert _{\mathcal{I%
}^{dual}\mathcal{\circ }LipL_{0}}.
\end{equation*}
\end{theorem}

\begin{proof}
Let $T\in \mathcal{I}_{dual}^{LipL}(X\times E;F)$. We aim to show that $%
T^{t}\in \mathcal{I}(F^{\ast };LipL_{0}(X\times E))$. Since the adjoint
operator $(\widehat{T})^{\ast }=Q^{-1}\circ T^{t}\in \mathcal{I}(F^{\ast };(X%
\widehat{\boxtimes }_{\pi }E)^{\ast }),$ where $Q:(X\widehat{\boxtimes }%
_{\pi }E)^{\ast }\rightarrow LipL_{0}(X\times E),$ is the canonical
isometric isomorphism given by%
\begin{equation*}
Q(\varphi )(x,e)=\varphi (\delta _{\left( x,0\right) }\boxtimes e),
\end{equation*}%
it follows that $\widehat{T}\in \mathcal{I}^{dual}(X\widehat{\boxtimes }%
_{\pi }E;F),$ and thus 
\begin{equation*}
T\in \mathcal{I}^{dual}\circ LipL_{0}(X\times E;F).
\end{equation*}%
Conversely, suppose $T\in \mathcal{I}^{dual}\circ LipL_{0}(X\times E;F)$.
Then there exist a Banach space $G,$ a Lip-linear operator $S\in
LipL_{0}(X\times E;G)$ and a linear operator $u$ such that $u^{\ast }\in 
\mathcal{I}(F^{\ast };G^{\ast })$ and 
\begin{equation*}
T=u\circ S.
\end{equation*}%
Taking duals, we have%
\begin{equation*}
T^{t}=S^{t}\circ u^{\ast }.
\end{equation*}%
By the ideal property, $T^{t}\in \mathcal{I}(F^{\ast };LipL_{0}(X\times E))$%
, which implies that $T\in \mathcal{I}_{dual}^{LipL}(X\times E;F)$. For the
relation $\left( \ref{2.1}\right) $, we have%
\begin{eqnarray*}
\Vert T\Vert _{\mathcal{I}_{dual}^{LipL}} &=&\Vert T^{t}\Vert _{\mathcal{I}%
}=\Vert Q\circ \left( \widehat{T}\right) ^{\ast }\Vert _{\mathcal{I}} \\
&\leq &\Vert \left( \widehat{T}\right) ^{\ast }\Vert _{\mathcal{I}}\Vert
Q\Vert =\Vert \widehat{T}\Vert _{\mathcal{I}^{dual}}=\Vert T\Vert _{\mathcal{%
I}^{dual}\mathcal{\circ }LipL_{0}}.
\end{eqnarray*}%
On the other hand, 
\begin{eqnarray*}
\Vert T\Vert _{\mathcal{I}^{dual}\mathcal{\circ }LipL_{0}} &=&\Vert \widehat{%
T}\Vert _{\mathcal{I}^{dual}}=\Vert \left( \widehat{T}\right) ^{\ast }\Vert
_{\mathcal{I}}=\Vert Q^{-1}\circ T^{t}\Vert _{\mathcal{I}} \\
&\leq &\Vert T^{t}\Vert _{\mathcal{I}}\Vert Q^{-1}\Vert =\Vert T\Vert _{%
\mathcal{I}_{dual}^{LipL}}.
\end{eqnarray*}%
This completes the argument.
\end{proof}

\begin{proposition}
Let $\mathit{{\mathcal{I}}}$ be a Banach operator ideal. Let $T\in
LipL_{0}(X\times E;F).$ The following assertions are equivalent:

1) The Lip-linear operator $T$ belongs to $\mathcal{I}\circ LipL_{0}(X\times
E;F).$

2) The bilinear operator $T_{bil}$ belongs to $\mathcal{I}\circ \mathcal{B}%
\left( \mathcal{F}\left( X\right) \times E;F\right) .$

Moreover, there is an isometric identification%
\begin{equation}
\mathcal{I}\circ LipL_{0}(X\times E;F)=\mathcal{I}\circ \mathcal{B}\left( 
\mathcal{F}\left( X\right) \times E;F\right) .  \label{2.4}
\end{equation}
\end{proposition}

\begin{proof}
The identification in $\left( \ref{1.4}\right) $ was established via the
correspondence $A\leftrightarrow A_{bil}.$ Using the relation $u\circ
A\leftrightarrow u\circ A_{bil},$ we can similarly establish the
identification in $\left( \ref{2.4}\right) .$ Moreover, by $\left( \ref{1.5}%
\right) $ $\widehat{T}=\widetilde{T_{bil}}\circ I.$ Then%
\begin{eqnarray*}
\Vert T\Vert _{\mathcal{I}\circ LipL_{0}} &=&\left\Vert \widehat{T}%
\right\Vert _{\mathcal{I}}=\left\Vert \widetilde{T_{bil}}\circ I\right\Vert
_{\mathcal{I}} \\
&\leq &\left\Vert \widetilde{T_{bil}}\right\Vert _{\mathcal{I}}\left\Vert
I\right\Vert =\left\Vert T_{bil}\right\Vert _{\mathcal{I}\circ \mathcal{B}}.
\end{eqnarray*}%
On the other hand, 
\begin{eqnarray*}
\left\Vert T_{bil}\right\Vert _{\mathcal{I}\circ \mathcal{B}} &=&\left\Vert 
\widetilde{T_{bil}}\right\Vert _{\mathcal{I}}=\left\Vert \widehat{T}\circ
I^{-1}\right\Vert _{\mathcal{I}} \\
&\leq &\left\Vert \widehat{T}\right\Vert _{\mathcal{I}}\left\Vert
I^{-1}\right\Vert =\Vert T\Vert _{\mathcal{I}\circ LipL_{0}}.
\end{eqnarray*}
\end{proof}

\begin{corollary}
Let $\mathcal{I},\mathcal{J}$ be operator ideals. Let $X$ be a pointed
metric space and let $E,F$ be Banach spaces. The following statements are
equivalent.

1) If $\mathcal{I}\circ LipL_{0}\left( X\times E;F\right) \subset \mathcal{J}%
\circ LipL_{0}\left( X\times E;F\right) ,$ then 
\begin{equation*}
\mathcal{I}\circ Lip_{0}\left( X;F\right) \subset \mathcal{J}\circ
Lip_{0}\left( X;F\right) \text{ and }\mathcal{I}\left( E;F\right) \subset 
\mathcal{J}\left( E;F\right) .
\end{equation*}

2) If $\mathcal{I}\circ LipL_{0}\left( X\times E;F\right) =LipL_{0}\left(
X\times E;F\right) ,$ then 
\begin{equation*}
\mathcal{I\circ }Lip_{0}\left( X;F\right) =Lip_{0}\left( X;F\right) \text{
and }\mathcal{I}\left( E;F\right) =\mathcal{L}\left( E;F\right) .
\end{equation*}
\end{corollary}

\begin{proof}
$1)$ By the identification $\left( \ref{2.4}\right) ,$ we deduce that $%
\mathcal{I}\circ \mathcal{B}\left( \mathcal{F}\left( X\right) \times
E;F\right) \subset \mathcal{J}\circ \mathcal{B}\left( \mathcal{F}\left(
X\right) \times E;F\right) .$ According to \cite[Lemma 3.4]{PelCompo}, we
have%
\begin{equation*}
\mathcal{I}\circ \mathcal{L}\left( \mathcal{F}\left( X\right) ;F\right)
\subset \mathcal{J}\circ \mathcal{L}\left( \mathcal{F}\left( X\right)
;F\right) \text{ and }\mathcal{I}\left( E;F\right) \subset \mathcal{J}\left(
E;F\right) .
\end{equation*}%
Using the identification $\left( \ref{1.0}\right) $, we obtain the desired
result.

$2)$ From $\left( \ref{1.4}\right) $ and $\left( \ref{2.4}\right) $, we have 
$\mathcal{I}\circ \mathcal{B}\left( \mathcal{F}\left( X\right) \times
E;F\right) =\mathcal{B}\left( \mathcal{F}\left( X\right) \times E;F\right) .$
Using \cite[Lemma 3.4]{PelCompo}, we get 
\begin{equation*}
\mathcal{I}\left( \mathcal{F}\left( X\right) ;F\right) =\mathcal{L}\left( 
\mathcal{F}\left( X\right) ;F\right) \text{ and  }\mathcal{I}\left(
E;F\right) =\mathcal{L}\left( E;F\right) ,
\end{equation*}%
Again, by the identification $\left( \ref{1.0}\right) $, we obtain the
desired result.
\end{proof}

As a direct consequence, we obtain the following result.

\begin{corollary}
Let $\mathit{{\mathcal{I}}}$ be an operator ideal. Let $F$ be a Banach
space. The following statements are equivalent.

1) $id_{F}\in \mathcal{I}\left( F;F\right) .$

2) $\mathcal{I}\circ LipL_{0}\left( X\times E;F\right) =LipL_{0}\left(
X\times E;F\right) $ for every pointed metric space $X$ and Banach space $E$.
\end{corollary}

\section{Cohen strongly summing Lip-Linear operators}

Let $\mathcal{D}_{p}$ denote the Banach operator ideal of strongly $p$%
-summing linear operators introduced by Cohen \cite{coh}. Using the
composition method, we construct a new Lip-Linear ideal based on the
operator ideal $\mathcal{D}_{p}$.

\begin{definition}
Let $1<p\leq \infty .$ Let $X$ be a pointed metric space, and let $E,F$ be
Banach spaces. A mapping $T\in LipL_{0}(X\times E;F)$ is called Cohen
strongly $p$-summing if: for each $\left( x_{i}\right) _{i=1}^{n},\left(
y_{i}\right) _{i=1}^{n}\subset X$, $\left( e_{i}\right) _{i=1}^{n}\subset E$
and $\left( z_{i}^{\ast }\right) _{i=1}^{n}\subset F^{\ast }$ we have 
\begin{equation}
\sum\limits_{i=1}^{n}\left\vert \left\langle
T(x_{i},e_{i})-T(y_{i},e_{i}),z_{i}^{\ast }\right\rangle \right\vert \leq
C(\sum\limits_{i=1}^{n}\left( d(x_{i},y_{i})\left\Vert e_{i}\right\Vert
\right) ^{p})^{\frac{1}{p}}\sup_{z^{\ast \ast }\in B_{F^{\ast \ast
}}}(\sum\limits_{i=1}^{n}\left\vert z^{\ast \ast }(z_{i}^{\ast
})\right\vert ^{p^{\ast }})^{\frac{1}{p^{\ast }}}.  \label{3.1}
\end{equation}%
We denote by $\mathcal{D}_{p}^{LipL}(X\times E;F)$ the space of all Cohen
strongly $p$-summing Lip-linear operators from $X\times E$ to $F$. For $T\in 
\mathcal{D}_{p}^{LipL}(X\times E;F)$, we define 
\begin{equation*}
d_{p}^{LipL}(T)=\inf \left\{ C:\text{ }C\text{ satisfies }\left( \ref{3.1}%
\right) \right\} .
\end{equation*}%
This defines a norm on $\mathcal{D}_{p}^{LipL}(X\times E;F),$ and it can be
shown that the pair $\left( \mathcal{D}_{p}^{LipL},d_{p}^{LipL}(\cdot
)\right) $ forms a Banach space.
\end{definition}

Let us provide an example of a Cohen strongly $p$-summing Lip-Linear
operator.

\begin{example}
Let $1<p\leq \infty .$\ Let $u:E\longrightarrow F$ be a strongly $p$-summing
linear operator and $f\in X^{\#}.$ The mapping $T:X\times E\longrightarrow
F,\ $defined by $T(x,e)=f(x)u(e),$ is a Cohen strongly $p$-summing
Lip-Linear operator with 
\begin{equation*}
d_{p}^{LipL}(T)\leq Lip(f)d_{p}\left( u\right) .
\end{equation*}%
Indeed, let $\left( x_{i}\right) _{i=1}^{n},\left( y_{i}\right)
_{i=1}^{n}\subset X$, $\left( e_{i}\right) _{i=1}^{n}\subset E$ and $\left(
z_{i}^{\ast }\right) _{i=1}^{n}\subset F^{\ast }$ we have%
\begin{eqnarray*}
\sum\limits_{i=1}^{n}\left\vert \left\langle
T(x_{i},e_{i})-T(y_{i},e_{i}),z_{i}^{\ast }\right\rangle \right\vert 
&=&\sum\limits_{i=1}^{n}\left\vert \left\langle
f(x_{i})u(e_{i})-f(y_{i})u(e_{i}),z_{i}^{\ast }\right\rangle \right\vert  \\
&=&\sum\limits_{i=1}^{n}\left\vert \left\langle u(\left(
f(x_{i})-f(y_{i})\right) e_{i}),z_{i}^{\ast }\right\rangle \right\vert .
\end{eqnarray*}%
Since $u$ is strongly $p$-summing, we get%
\begin{eqnarray*}
&&\sum\limits_{i=1}^{n}\left\vert \left\langle u(\left(
f(x_{i})-f(y_{i})\right) e_{i}),z_{i}^{\ast }\right\rangle \right\vert  \\
&\leq &d_{p}\left( u\right) (\sum\limits_{i=1}^{n}\left( \left\Vert \left(
f(x_{i})-f(y_{i})\right) e_{i}\right\Vert \right) ^{p})^{\frac{1}{p}%
}\sup_{z^{\ast \ast }\in B_{F^{\ast \ast
}}}(\sum\limits_{i=1}^{n}\left\vert z^{\ast \ast }(z_{i}^{\ast
})\right\vert ^{p^{\ast }})^{\frac{1}{p^{\ast }}} \\
&\leq &Lip\left( f\right) d_{p}\left( u\right) (\sum\limits_{i=1}^{n}\left(
d\left( x_{i},y_{i}\right) \left\Vert e_{i}\right\Vert \right) ^{p})^{\frac{1%
}{p}}\sup_{z^{\ast \ast }\in B_{F^{\ast \ast
}}}(\sum\limits_{i=1}^{n}\left\vert z^{\ast \ast }(z_{i}^{\ast
})\right\vert ^{p^{\ast }})^{\frac{1}{p^{\ast }}}.
\end{eqnarray*}
\end{example}

\begin{theorem}
\label{TheoremDomination}Let $1<p\leq \infty .$ Consider $T\in
LipL_{0}(X\times E;F)$. The following assertions are equivalent.\newline
1) The Lip-Linear operator $T$ is Cohen strongly $p$-summing.\newline
2) There exist a positive constant $C$ and a Borel probability measure $\mu $
on $B_{F^{\ast \ast }}$ such that for all $x,y\in X$, $e\in E$ and $z^{\ast
}\in F^{\ast }$ 
\begin{equation}
\left\vert \left\langle T\left( x,e\right) -T\left( y,e\right) ,z^{\ast
}\right\rangle \right\vert \leq Cd(x,y)\left\Vert e\right\Vert \left(
\int_{B_{F^{\ast \ast }}}|\langle z^{\ast \ast },z^{\ast }\rangle |^{p^{\ast
}}d\mu \right) ^{\frac{1}{p^{\ast }}}.  \label{3.1.2}
\end{equation}%
Moreover, in this case $d_{p}^{LipL}(T)=\inf \left\{ C:\text{ }C\text{
satisfies }\left( \ref{3.1.2}\right) \right\} $.
\end{theorem}

\begin{proof}
Define the following mappings:%
\begin{equation*}
\left\{ 
\begin{array}{l}
R_{1}:B_{F^{\ast \ast }}\times (X\times X\times E)\times \mathbb{R}%
\rightarrow \left[ 0;\infty \right[ :R_{1}(z^{\ast \ast },(x,y,e),\lambda
)=\left\vert \lambda \right\vert d(x,y)\Vert e\Vert  \\ 
R_{2}:B_{F^{\ast \ast }}\times (X\times X\times E)\times F^{\ast
}\rightarrow \left[ 0;\infty \right[ :R_{2}(z^{\ast \ast },(x,y,e),z^{\ast
})=|z^{\ast \ast }(z^{\ast })| \\ 
S:LipL_{0}(X\times E;F)\times (X\times X\times E)\times \mathbb{R}\times
F^{\ast }\rightarrow \left[ 0;\infty \right[ : \\ 
S(T,(x,y,e),\lambda ,z^{\ast })=\left\vert \lambda \right\vert \left\vert
\left\langle T(x,e)-T(y,e),z^{\ast }\right\rangle \right\vert .%
\end{array}%
\right. 
\end{equation*}%
These mappings satisfy conditions $\left( 1\right) $ and $\left( 2\right) $
in \cite[Page 1255]{pss}. Therefore, $T:X\times E\rightarrow F$ is Cohen
strongly $p$-summing if, and only if, 
\begin{eqnarray*}
&&\sum\limits_{i=1}^{n}S(T,(x_{i},y_{i},e_{i}),\lambda _{i},z_{i}^{\ast }) \\
&\leq &C(\sum\limits_{i=1}^{n}R_{1}(z^{\ast \ast
},(x_{i},y_{i},e_{i}),\lambda _{i})^{p})^{\frac{1}{p}}\sup\limits_{z^{\ast
\ast }\in B_{E^{\ast \ast }}}(\sum\limits_{i=1}^{n}R_{2}(z^{\ast \ast
},(x_{i},y_{i},e_{i}),z_{i}^{\ast })^{p^{\ast }})^{\frac{1}{p^{\ast }}}.
\end{eqnarray*}%
Thus, $T$ is $R_{1},R_{2}$-$S$-abstract $(p,p^{\ast })$-summing. By \cite[%
Theorem 4.6]{pss}, this is equivalent to the existence of a constant $C$ and
probability measures $\mu _{1}$ and $\mu _{2}$ on $B_{F^{\ast \ast }}$ such
that%
\begin{eqnarray*}
&&S(T,(x,y,e),\lambda ,z^{\ast }) \\
&\leq &C\left( \int\nolimits_{B_{F^{\ast \ast }}}R_{1}(z^{\ast \ast
},(x,y,e),\lambda )^{p}d\mu _{1}\right) ^{\frac{1}{p}}\left(
\int\nolimits_{B_{F^{\ast \ast }}}R_{2}(z^{\ast \ast },(x,y),z^{\ast
})^{p^{\ast }}d\mu _{2}\right) ^{\frac{1}{p^{\ast }}}.
\end{eqnarray*}%
Consequently, 
\begin{equation*}
\left\vert \left\langle T(x,e)-T\left( y,e\right) ,z^{\ast }\right\rangle
\right\vert \leq Cd(x,y)\left\Vert e\right\Vert \left( \int_{B_{F^{\ast \ast
}}}|\langle z^{\ast \ast },z^{\ast }\rangle |^{p^{\ast }}d\mu \right) ^{%
\frac{1}{p^{\ast }}},
\end{equation*}%
where $\mu :=\mu _{2}$.
\end{proof}

The multilinear definition of strongly $p$-summing operators was introduced
by Achour and Mezrag in \cite{achmez}. In the result below, we show that
this definition coincides with the one for Lip-linear operators when
restricted to the bilinear case.

\begin{proposition}
Let $1<p\leq \infty .$\ Let $E,F$ and $G$ be Banach spaces. Let $T:E\times
F\rightarrow G$ be a bilinear operator. Then, $T\in \mathcal{D}%
_{p}^{LipL}(E\times F;G)$ if, and only if, $T\in \mathcal{D}_{p}^{2}(E\times
F;G).$ Moreover, we have 
\begin{equation*}
d_{p}^{2}(T)=d_{p}^{LipL}(T).
\end{equation*}
\end{proposition}

\begin{proof}
Suppose $T\in \mathcal{D}_{p}(E\times F;G).$ For $x,y\in E$ and $e\in F,$%
\begin{equation*}
\left\vert \left\langle T(x,e)-T(y,e\rangle ,z^{\ast }\right\rangle
\right\vert =\left\vert \left\langle T(x-y,e),z^{\ast }\right\rangle
\right\vert 
\end{equation*}%
Applying the integral characterization of the class $\mathcal{D}_{p}^{2}$
given in \cite[Theorem 2.4]{achmez}, we obtain%
\begin{equation*}
\left\vert \left\langle T(x,e)-T(y,e\rangle ,z^{\ast }\right\rangle
\right\vert \leq d_{p}^{2}(T)\left\Vert x-y\right\Vert \left\Vert
e\right\Vert \left( \int_{B_{F^{\ast \ast }}}|\langle z^{\ast \ast },z^{\ast
}\rangle |^{p^{\ast }}d\mu \right) ^{\frac{1}{p^{\ast }}}.
\end{equation*}%
This shows that $T\in \mathcal{D}_{p}^{LipL}(E\times F;G)$ and $%
d_{p}^{LipL}(T)\leq d_{p}^{2}(T).$ Conversely, suppose $T\in \mathcal{D}%
_{p}^{LipL}(E\times F;G).$ For $x\in E,e\in F$ and $z^{\ast }\in G^{\ast }$ 
\begin{equation*}
\left\vert \left\langle T(x,e),z^{\ast }\right\rangle \right\vert
=\left\vert \left\langle T(x,e)-T(0,e\rangle ,z^{\ast }\right\rangle
\right\vert .
\end{equation*}%
By Theorem $\ref{TheoremDomination},$\ we have%
\begin{equation*}
\left\vert \left\langle T(x,e),z^{\ast }\right\rangle \right\vert \leq
d_{p}^{LipL}(T)\left\Vert x-0\right\Vert \Vert e\Vert \left(
\int_{B_{F^{\ast \ast }}}|\langle z^{\ast \ast },z^{\ast }\rangle |^{p^{\ast
}}d\mu \right) ^{\frac{1}{p^{\ast }}}.
\end{equation*}%
Thus, $T\in \mathcal{D}_{p}^{2}(E\times F;G)$ and $d_{p}^{2}(T)\leq
d_{p}^{LipL}(T).$
\end{proof}

Let $T:X\times E\longrightarrow F$ be a Lip-linear operator. We define the
associated operators:%
\begin{eqnarray*}
A_{T} &:&x\mapsto A_{T}\left( x\right) ,\text{ where }A_{T}(x)\left(
e\right) =T(x,e)\text{ for every }x\in X\text{ and }e\in E. \\
B_{T} &:&e\mapsto B_{T}\left( e\right) ,\text{ where }B_{T}(e)\left(
x\right) =T(x,e)\text{ for every }x\in X\text{ and }e\in E.
\end{eqnarray*}

A useful characterization linking a Lip-linear operator with its associated
mappings is provided in the following corollary.

\begin{proposition}
\label{Propo2}Let $1<p\leq \infty $ and $T\in LipL_{0}(X\times E;F)$, the
following assertions are equivalent.

1) The Lip-linear operator $T$ belongs to $\mathcal{D}_{p}^{LipL}\left(
X\times E;F\right) $.

2) The Lipschitz operator $A_{T}$ belongs to $Lip_{0}(X;\mathcal{D}%
_{p}(E;F)) $.

3) The linear operator $B_{T}$ belongs to $\mathcal{L}(E;\mathcal{D}%
_{p}^{Lip}(X;F))$.

Additionally, the following equality hold%
\begin{equation*}
d_{p}^{LipL}(T)=Lip(A_{T})=\left\Vert B_{T}\right\Vert .
\end{equation*}%
As a consequence, we obtain the isometric identifications%
\begin{equation*}
\mathcal{D}_{p}^{LipL}\left( X\times E;F\right) =Lip_{0}(X;\mathcal{D}%
_{p}(E;F))=\mathcal{L}(E;\mathcal{D}_{p}^{Lip}(X;F)).
\end{equation*}
\end{proposition}

\begin{proof}
$1)\Longrightarrow 2):$ Let $T\in \mathcal{D}_{p}^{LipL}\left( X\times
E;F\right) $ and fix $x,y\in X.$ We aim to show that $A_{T}(x)-A_{T}\left(
y\right) \in \mathcal{D}_{p}(E;F).$ Let $e\in E$ and $z^{\ast }\in F^{\ast
}, $ we have 
\begin{equation*}
\left\vert \left\langle \left( A_{T}(x)-A_{T}\left( y\right) \right) \left(
e\right) ,z^{\ast }\right\rangle \right\vert =\left\vert \left\langle
T(x,e)-T(y,e),z^{\ast }\right\rangle \right\vert .
\end{equation*}%
Since $T$ is Cohen strongly $p$-summing, we obtain%
\begin{eqnarray*}
&&\left\vert \left\langle \left( A_{T}(x)-A_{T}\left( y\right) \right)
\left( e\right) ,z^{\ast }\right\rangle \right\vert \\
&\leq &d_{p}^{LipL}(T)d\left( x,y\right) \left\Vert e\right\Vert
(\int\limits_{B_{F^{\ast \ast }}}\left\vert z^{\ast \ast }\left( z^{\ast
}\right) \right\vert ^{p^{\ast }}d\mu )^{\frac{1}{p^{\ast }}}
\end{eqnarray*}%
This shows that $A_{T}(x)-A_{T}\left( y\right) \in \mathcal{D}_{p}(E;F)$ with%
\begin{equation}
d_{p}\left( A_{T}(x)-A_{T}\left( y\right) \right) \leq
d_{p}^{LipL}(T)d\left( x,y\right) .  \label{3.6}
\end{equation}%
Then, $A_{T}$ is Lipschitz and 
\begin{equation*}
Lip\left( A_{T}\right) \leq d_{p}^{LipL}\left( T\right) .
\end{equation*}%
Setting $y=0$ in $\left( \ref{3.6}\right) $ implies $A_{T}(x)\in \mathcal{D}%
_{p}(E;F)$.

$2)\Longrightarrow 3):\mathcal{\ }$Let $e\in E.$ We aim to show $B_{T}\left(
e\right) \in \mathcal{D}_{p}^{Lip}(X;F).$ For $x,y\in X$ and $z^{\ast }\in
F^{\ast }.$ Using the Pietch Domination Theorem of strongly $p$-summing
operator of $A_{T}(x)-A_{T}\left( y\right) $ (see \cite[Theorem 2.3.1]{coh})$%
,$ there exists a Borel probability measure $\mu $ on $B_{F^{\ast \ast }}$
such that%
\begin{eqnarray*}
\left\vert \left\langle B_{T}\left( e\right) \left( x\right) -B_{T}\left(
e\right) \left( y\right) ,z^{\ast }\right\rangle \right\vert  &=&\left\vert
\left\langle \left( A_{T}(x)-A_{T}\left( y\right) \right) \left( e\right)
,z^{\ast }\right\rangle \right\vert  \\
&\leq &d_{p}\left( A_{T}(x)-A_{T}\left( y\right) \right) \left\Vert
e\right\Vert (\int\limits_{B_{F^{\ast \ast }}}\left\vert z^{\ast \ast
}\left( z^{\ast }\right) \right\vert ^{p^{\ast }}d\mu )^{\frac{1}{p^{\ast }}}
\\
&\leq &Lip\left( A_{T}\right) d\left( x,y\right) \left\Vert e\right\Vert
(\int\limits_{B_{F^{\ast \ast }}}\left\vert z^{\ast \ast }\left( z^{\ast
}\right) \right\vert ^{p^{\ast }}d\mu )^{\frac{1}{p^{\ast }}}.
\end{eqnarray*}%
This implies that $B_{T}\left( e\right) $ is Lipschitz strongly $p$-summing,
and%
\begin{equation*}
d_{p}^{L}\left( B_{T}\left( e\right) \right) \leq Lip\left( A_{T}\right)
\left\Vert e\right\Vert .
\end{equation*}%
Thus, $B_{T}$ is bounded, and%
\begin{equation*}
\left\Vert B_{T}\right\Vert \leq Lip\left( A_{T}\right) .
\end{equation*}%
$3)\Longrightarrow 1):$ Let $x,y\in X$ and $e\in E$. Then,%
\begin{equation*}
\left\vert \left\langle T(x,e)-T(y,e),z^{\ast }\right\rangle \right\vert
=\left\vert \left\langle B_{T}\left( e\right) \left( x\right) -B_{T}\left(
e\right) \left( y\right) ,z^{\ast }\right\rangle \right\vert .
\end{equation*}%
Since $B_{T}\left( e\right) $ is Lipschitz strongly $p$-summing, there
exists a Borel probability measure $\mu $ on $B_{F^{\ast \ast }}$ such that%
\begin{eqnarray*}
\left\vert \left\langle T(x,e)-T(y,e),z^{\ast }\right\rangle \right\vert 
&=&\left\vert \left\langle B_{T}\left( e\right) \left( x\right) -B_{T}\left(
e\right) \left( y\right) ,z^{\ast }\right\rangle \right\vert  \\
&\leq &d_{p}^{L}(B_{T}\left( e\right) )d\left( x,y\right)
(\int\limits_{B_{F^{\ast \ast }}}\left\vert z^{\ast \ast }\left( z^{\ast
}\right) \right\vert ^{p^{\ast }}d\mu )^{\frac{1}{p^{\ast }}} \\
&\leq &\left\Vert B_{T}\right\Vert d\left( x,y\right) \left\Vert
e\right\Vert (\int\limits_{B_{F^{\ast \ast }}}\left\vert z^{\ast \ast
}\left( z^{\ast }\right) \right\vert ^{p^{\ast }}d\mu )^{\frac{1}{p^{\ast }}%
}.
\end{eqnarray*}%
Thus, $T$ is Cohen strongly $p$-summing, and%
\begin{equation*}
d_{p}^{LipL}(T)\leq \left\Vert B_{T}\right\Vert .
\end{equation*}%
We show the surjectivity. Let $S\in Lip_{0}(X;\mathcal{D}_{p}(E;F))$. Define
a Lip-linear operator $T_{S}:X\times E\rightarrow F$ by%
\begin{equation*}
T_{S}\left( x,e\right) =S\left( x\right) \left( e\right) .
\end{equation*}%
We claim that $T_{S}$ is Cohen strongly $p$-summing. Indeed, let $x,y\in
X,e\in E$ and $z^{\ast }\in F^{\ast }$. Then,%
\begin{eqnarray*}
\left\vert \left\langle T_{S}(x,e)-T_{S}(y,e),z^{\ast }\right\rangle
\right\vert  &=&\left\vert \left\langle S(x)\left( e\right) -S(y)\left(
e\right) ,z^{\ast }\right\rangle \right\vert  \\
&\leq &d_{p}\left( S(x)-S(y)\right) \left\Vert e\right\Vert
(\int\limits_{B_{F^{\ast \ast }}}\left\vert z^{\ast \ast }\left( z^{\ast
}\right) \right\vert ^{p^{\ast }}d\mu )^{\frac{1}{p^{\ast }}} \\
&\leq &Lip\left( S\right) d\left( x,y\right) \left\Vert e\right\Vert
(\int\limits_{B_{F^{\ast \ast }}}\left\vert z^{\ast \ast }\left( z^{\ast
}\right) \right\vert ^{p^{\ast }}d\mu )^{\frac{1}{p^{\ast }}}.
\end{eqnarray*}%
This shows that $T_{S}$ is strongly $p$-summing, and therefore $S$ induces a
surjective correspondence. Using similar arguments, a corresponding
identification holds for operators in $\mathcal{L}(E;\mathcal{D}%
_{p}^{L}(X;F)).$
\end{proof}

\begin{corollary}
Let $1<p\leq \infty .$ Consider $T\in LipL_{0}(X\times E;F)$. The following
assertions are equivalent.

1) The Lip-Linear operator $T$ is Cohen strongly $p$-summing.

2) There exist a Banach space $G$ and a $p^{\ast }$-summing linear operator $%
v:F^{\ast }\longrightarrow G$ such that, for all $x,y\in X,$\ $e\in E$ and $%
z^{\ast }\in F^{\ast },$ the following inequality holds 
\begin{equation}
\left\vert \left\langle T\left( x,e\right) -T\left( y,e\right) ,z^{\ast
}\right\rangle \right\vert \leq d(x,y)\left\Vert e\right\Vert \left\Vert
v(z^{\ast })\right\Vert .  \label{3.2}
\end{equation}

In this case, we have%
\begin{equation*}
d_{p}^{LipL}(T)=\inf \left\{ \pi _{p^{\ast }}\left( v\right) \right\} ,
\end{equation*}%
where the infimum is taken over all Banach space $G$ and operator $v$ such
that the inequality $\left( \ref{3.2}\right) $ is satisfied.
\end{corollary}

\begin{proof}
$1)\Longrightarrow 2):$ Let $T\in LipL_{0}(X\times E;F)$ be a Cohen strongly 
$p$-summing. Let $i_{F^{\ast }}$ be the natural isometric embedding of $%
F^{\ast }$ in $\mathcal{C}(B_{F^{\ast \ast }}),$ composed with the formal
identity map from $\mathcal{C}(B_{F^{\ast \ast }})$ into $L_{\infty }(\mu ),$
defined as 
\begin{equation*}
i_{F^{\ast }}(z^{\ast })(z^{\ast \ast })=z^{\ast \ast }(z^{\ast }),
\end{equation*}%
for $z^{\ast }\in F^{\ast },\ z^{\ast \ast }\in {F^{\ast \ast }}$. Let $%
I_{\infty ,p^{\ast }}:L_{\infty }(\mu )\longrightarrow L_{p^{\ast }}(\mu )$
be the canonical mapping $I_{\infty ,p^{\ast }}(f)=f$. Note that (see \cite[%
Examples 2.9]{distel}), $I_{\infty ,p^{\ast }}$ is $p^{\ast }$-summing and $%
\pi _{p^{\ast }}(I_{\infty ,p^{\ast }})=1$. By Theorem \ref%
{TheoremDomination}, we have%
\begin{equation*}
\left\vert \left\langle T(x,e)-T\left( y,e\right) ,z^{\ast }\right\rangle
\right\vert \leq d_{p}^{LipL}(T)d(x,y)\left\Vert e\right\Vert \left(
\int_{B_{F^{\ast \ast }}}|\left\langle z^{\ast \ast },z^{\ast }\right\rangle
|^{p^{\ast }}d\mu \right) ^{\frac{1}{p^{\ast }}}.
\end{equation*}%
This inequality can be rewritten as%
\begin{equation*}
\left\vert \left\langle T(x,e)-T\left( y,e\right) ,z^{\ast }\right\rangle
\right\vert \leq d(x,y)\left\Vert e\right\Vert \left( \int_{B_{F^{\ast \ast
}}}\left\vert d_{p}^{LipL}(T)I_{\infty ,p^{\ast }}i_{F^{\ast }}\left(
z^{\ast }\right) \left( z^{\ast \ast }\right) \right\vert ^{p^{\ast }}d\mu
\right) ^{\frac{1}{p^{\ast }}}.
\end{equation*}%
Thus, there exists a Banach space $F=L_{p^{\ast }}(\mu )$ and a $p^{\ast }$%
-summing operator $v=d_{p}^{LipL}(T)I_{\infty ,p^{\ast }}i_{F^{\ast }},$
such that 
\begin{equation*}
\left\vert \left\langle T(x,e)-T\left( y,e\right) ,z^{\ast }\right\rangle
\right\vert \leq d\left( x,y\right) \left\Vert e\right\Vert \left\Vert
v\left( z^{\ast }\right) \right\Vert .
\end{equation*}%
We have%
\begin{eqnarray*}
\pi _{p^{\ast }}\left( v\right)  &=&\pi _{p^{\ast }}\left(
d_{p}^{LipL}(T)I_{\infty ,p^{\ast }}i_{F^{\ast }}\right)  \\
&\leq &d_{p}^{LipL}(T)\pi _{p^{\ast }}(I_{\infty ,p^{\ast }})\left\Vert
i_{F^{\ast }}\right\Vert =d_{p}^{LipL}(T).
\end{eqnarray*}%
Therefore,%
\begin{equation*}
\inf \left\{ \pi _{p^{\ast }}\left( v\right) \right\} \leq d_{p}^{LipL}(T)
\end{equation*}%
$2)\Longrightarrow 1):$ Assume there exist a Banach space $G$ and$\ v\in \Pi
_{p}(F^{\ast };G)$ such that 
\begin{equation*}
\left\vert \left\langle T(x,e)-T\left( y,e\right) ,z^{\ast }\right\rangle
\right\vert \leq d(x,y)\left\Vert e\right\Vert \left\Vert v\left( z^{\ast
}\right) \right\Vert ,
\end{equation*}%
for all $x,y\in X,$ $e\in E$ and $z^{\ast }\in F^{\ast }.$ Since $v$ is $%
p^{\ast }$-summing, there exists a regular Borel probability measure $\mu $
on $B_{F^{\ast \ast }}$ such that 
\begin{equation*}
\left\vert \left\langle T(x,e)-T\left( y,e\right) ,z^{\ast }\right\rangle
\right\vert \leq \pi _{p^{\ast }}\left( v\right) d(x,y)\left\Vert
e\right\Vert \left( \int_{B_{F^{\ast \ast }}}|\langle z^{\ast \ast },z^{\ast
}\rangle |^{p^{\ast }}d\mu \right) ^{\frac{1}{p^{\ast }}}.
\end{equation*}%
By Theorem \ref{TheoremDomination}, $T$ is Cohen strongly $p$-summing, with%
\begin{equation*}
d_{p}^{LipL}(T)\leq \pi _{p^{\ast }}\left( v\right) .
\end{equation*}%
This holds for every $v$ that satisfies inequality $\left( \ref{3.2}\right) $%
. Therefore,%
\begin{equation*}
d_{p}^{LipL}(T)\leq \inf \left\{ \pi _{p^{\ast }}\left( v\right) \right\} .
\end{equation*}
\end{proof}

\begin{theorem}
Let $X$ be a pointed metric space, and let $E,F$ be Banach spaces. For $%
1<p\leq \infty ,$ the following assertions are equivalent.

1) The lip-Linear operator $T\in \mathcal{D}_{p}^{LipL}(X\times E;F)$

2) The dual operator $T^{t}\in \Pi _{p^{\ast }}(F^{\ast };LipL_{0}\left(
X\times E\right) ).$

In this case, we have%
\begin{equation*}
d_{p}^{LipL}(T)=\pi _{p^{\ast }}(T^{t}).
\end{equation*}%
\end{theorem}

\begin{proof}
$1)\Longrightarrow 2):$ Assume that $T$ is\ Cohen strongly $p$-summing. Let $%
z^{\ast }\in F^{\ast }.$ The norm of $T^{t}\left( z^{\ast }\right) $
satisfies%
\begin{eqnarray*}
LipL\left( T^{t}\left( z^{\ast }\right) \right)  &=&\sup_{x\neq y,\left\Vert
e\right\Vert =1}\frac{\left\vert T^{t}\left( z^{\ast }\right) \left(
x,e\right) -T^{t}\left( z^{\ast }\right) \left( y,e\right) \right\vert }{%
d\left( x,y\right) } \\
&=&\sup_{x\neq y,\left\Vert e\right\Vert =1}\frac{\left\vert \left\langle
T\left( x,e\right) -T\left( y,e\right) ,z^{\ast }\right\rangle \right\vert }{%
d\left( x,y\right) }
\end{eqnarray*}%
Since $T$ is Cohen strongly $p$-summing, by Theorem $\ref{TheoremDomination},
$ there exists a Borel probability measure $\mu $ on $B_{F^{\ast \ast }}$
such that%
\begin{equation*}
\left\vert \left\langle T\left( x,e\right) -T\left( y,e\right) ,z^{\ast
}\right\rangle \right\vert \leq d_{p}^{LipL}\left( T\right) d(x,y)\left\Vert
e\right\Vert \left( \int_{B_{F^{\ast \ast }}}|\langle z^{\ast \ast },z^{\ast
}\rangle |^{p^{\ast }}d\mu \right) ^{\frac{1}{p^{\ast }}}.
\end{equation*}%
Thus,%
\begin{equation*}
LipL\left( T^{t}\left( z^{\ast }\right) \right) \leq d_{p}^{LipL}\left(
T\right) \left( \int_{B_{F^{\ast \ast }}}|\langle z^{\ast \ast },z^{\ast
}\rangle |^{p^{\ast }}d\mu \right) ^{\frac{1}{p^{\ast }}}
\end{equation*}%
This implies that $T^{t}$ is $p^{\ast }$-summing, and 
\begin{equation*}
\pi _{p^{\ast }}\left( T^{t}\right) \leq d_{p}^{LipL}\left( T\right) .
\end{equation*}%
$2)\Longrightarrow 1):$ Suppose that $T^{t}$ is $p^{\ast }$-summing. Then
for any $x,y\in X,e\in E$ and $z^{\ast }\in F^{\ast },$ we have%
\begin{eqnarray*}
\left\vert \left\langle T\left( x,e\right) -T\left( y,e\right) ,z^{\ast
}\right\rangle \right\vert  &=&\left\vert T^{t}\left( z^{\ast }\right)
\left( x,e\right) -T^{t}\left( z^{\ast }\right) \left( y,e\right)
\right\vert  \\
&\leq &LipL\left( T^{t}\left( z^{\ast }\right) \right) d\left( x,y\right)
\left\Vert e\right\Vert .
\end{eqnarray*}%
Using the fact that $T^{t}\in \Pi _{p^{\ast }}(F^{\ast };LipL_{0}\left(
X\times E\right) ),$ there exists a Radon probability measure $\mu $ on $%
B_{F^{\ast \ast }}$ such that%
\begin{equation*}
LipL\left( T^{t}\left( z^{\ast }\right) \right) \leq \pi _{p^{\ast }}\left(
T^{t}\right) \left( \int_{B_{F^{\ast \ast }}}|\langle z^{\ast \ast },z^{\ast
}\rangle |^{p^{\ast }}d\mu \right) ^{\frac{1}{p^{\ast }}}
\end{equation*}%
Hence, 
\begin{equation*}
\left\vert \left\langle T\left( x,e\right) -T\left( y,e\right) ,z^{\ast
}\right\rangle \right\vert \leq \pi _{p^{\ast }}\left( T^{t}\right) d\left(
x,y\right) \left\Vert e\right\Vert \left( \int_{B_{F^{\ast \ast }}}|\langle
z^{\ast \ast },z^{\ast }\rangle |^{p^{\ast }}d\mu \right) ^{\frac{1}{p^{\ast
}}}
\end{equation*}%
This shows that $T$ is Cohen strongly $p$-summing and%
\begin{equation*}
d_{p}^{LipL}\left( T\right) \leq \pi _{p^{\ast }}\left( T^{t}\right) .
\end{equation*}
\end{proof}

As a consequence of the above Theorem and Proposition $\ref{Propo1}$, the
class $\mathcal{D}_{p}^{LipL}$ constitutes a Banach ideal of Lip-linear
operators. Moreover, for every pointed metric space $X$ and Banach spaces $%
E,F$ we have the isometric identification%
\begin{equation*}
\mathcal{D}_{p}^{LipL}(X\times E;F)=\left( \Pi _{p^{\ast }}\right)
_{dual}^{LipL}(X\times E;F).
\end{equation*}

We now show that $\mathcal{D}_{p}^{LipL}$ arises naturally by applying the
composition method to the classical operator ideal $\mathcal{D}_{p}$. To
establish this, we analyze the relationship between a Lip-Linear operator
and its linearization for the concept of strongly $p$-summing operators.

\begin{theorem}
\label{T_L}Let $X$ be a pointed metric space, and let $E,F$ be Banach
spaces. For $1<p\leq \infty ,$ we have $T\in \mathcal{D}_{p}^{LipL}(X\times
E;F)$ if, and only if, its linearization $\widehat{T}\in \mathcal{D}_{p}(X%
\widehat{\boxtimes }_{\pi }E;F).$ In this case,%
\begin{equation*}
d_{p}^{LipL}(T)=d_{p}(\widehat{T}).
\end{equation*}
\end{theorem}

\begin{proof}
Suppose that $\widehat{T}$ is strongly $p$-summing. By \cite[Theorem 2.2.2]%
{coh}$,$ the dual operator $\left( \widehat{T}\right) ^{\ast }:F^{\ast
}\rightarrow \left( X\widehat{\boxtimes }_{\pi }E\right) ^{\ast }$ is $%
p^{\ast }$-summing. To establish the relation $\left( \ref{3.2}\right) $, we
set $G=\left( X\widehat{\boxtimes }_{\pi }E\right) ^{\ast }$ and $v=\left( 
\widehat{T}\right) ^{\ast }.$ Let $x,y\in X$ and $e\in E,$ we have%
\begin{eqnarray*}
\left\vert \left\langle T\left( x,e\right) -T\left( y,e\right) ,z^{\ast
}\right\rangle \right\vert &=&\left\vert \left\langle \widehat{T}\left(
\delta _{x}\boxtimes e\right) -\widehat{T}\left( \delta _{y}\boxtimes
e\right) ,z^{\ast }\right\rangle \right\vert \\
&=&\left\vert \left\langle \widehat{T}\left( \delta _{\left( x,y\right)
}\boxtimes e\right) ,z^{\ast }\right\rangle \right\vert =\left\vert
\left\langle \delta _{\left( x,y\right) }\boxtimes e,\left( \widehat{T}%
\right) ^{\ast }\left( z^{\ast }\right) \right\rangle \right\vert \\
&\leq &\pi \left( \delta _{\left( x,y\right) }\boxtimes e\right) \left\Vert
\left( \widehat{T}\right) ^{\ast }(z^{\ast })\right\Vert \\
&\leq &d(x,y)\Vert e\Vert \left\Vert \left( \widehat{T}\right) ^{\ast
}(z^{\ast })\right\Vert .
\end{eqnarray*}%
Thus, $T\in \mathcal{D}_{p}^{LipL}(X\times E;F),$ and we have%
\begin{equation*}
d_{p}^{LipL}(T)\leq \pi _{p^{\ast }}(\left( \widehat{T}\right) ^{\ast
})=d_{p}(\widehat{T}).
\end{equation*}%
For the converse, suppose that $T\in \mathcal{D}_{p}^{LipL}(X\times E;F).$
Then, there exist a Banach space $G$ and a $p^{\ast }$-summing linear
operator $v:F^{\ast }\longrightarrow G$ such that the inequality $\left( \ref%
{3.2}\right) $ holds. Let $u\in X\widehat{\boxtimes }_{\pi }E$ and $z^{\ast
}\in F^{\ast }$. For $\varepsilon >0,$ we can choose a representation $%
u=\sum_{i=1}^{n}\delta _{\left( x_{i},y_{i}\right) }\boxtimes e_{i}$ such
that 
\begin{equation*}
\sum_{i=1}^{n}d(x_{i},y_{i})\Vert e_{i}\Vert \leq \pi (u)+\varepsilon .
\end{equation*}%
Then,%
\begin{eqnarray*}
\left\vert \left\langle \widehat{T}\left( u\right) ,z^{\ast }\right\rangle
\right\vert &=&\left\vert \sum_{i=1}^{n}\left\langle T\left(
x_{i},e_{i}\right) -T\left( y_{i},e_{i}\right) ,z^{\ast }\right\rangle
\right\vert \\
&\leq &\sum_{i=1}^{n}\left\vert \left\langle T\left( x_{i},e_{i}\right)
-T\left( y_{i},e_{i}\right) ,z^{\ast }\right\rangle \right\vert \\
&\leq &\sum_{i=1}^{n}d(x_{i},y_{i})\Vert e_{i}\Vert \left\Vert v(z^{\ast
})\right\Vert \\
&\leq &\left( \pi (u)+\varepsilon \right) \left\Vert v(z^{\ast })\right\Vert
.
\end{eqnarray*}%
By taking $\varepsilon \rightarrow 0$, we find 
\begin{equation*}
\left\vert \left\langle \widehat{T}\left( u\right) ,z^{\ast }\right\rangle
\right\vert \leq \pi (u)\left\Vert v(z^{\ast })\right\Vert .
\end{equation*}%
Since $v$ is $p^{\ast }$-summing, by \cite[Theorem 2.12]{distel}, there
exist a regular Borel probability measure $\mu $ on $B_{F^{\ast \ast }}$
such that 
\begin{equation*}
\left\Vert v(z^{\ast })\right\Vert \leq \pi _{p}\left( v\right) \left(
\int_{B_{F^{\ast \ast }}}|\langle z^{\ast \ast },z^{\ast }\rangle |^{p^{\ast
}}d\mu \right) ^{\frac{1}{p^{\ast }}}.
\end{equation*}%
Then, 
\begin{equation*}
\left\vert \left\langle \widehat{T}\left( u\right) ,z^{\ast }\right\rangle
\right\vert \leq \pi _{p^{\ast }}\left( v\right) \pi (u)\left(
\int_{B_{F^{\ast \ast }}}|\langle z^{\ast \ast },z^{\ast }\rangle |^{p^{\ast
}}d\mu \right) ^{\frac{1}{p^{\ast }}}.
\end{equation*}%
By the Pietsch Domination Theorem for strongly $p$-summing operators \cite[%
Theorem 2.3.1]{coh}, we conclude that $\widehat{T}$ is strongly $p$-summing,
and we obtain%
\begin{equation}
d_{p}(\widehat{T})\leq \pi _{p^{\ast }}\left( v\right) .  \label{3.3}
\end{equation}%
The inequality $\left( \ref{3.3}\right) $\ holds for every $v$ that
satisfies inequality $\left( \ref{3.2}\right) $. Therefore,%
\begin{equation*}
d_{p}(\widehat{T})\leq d_{p}^{Lipl}\left( T\right) .
\end{equation*}
\end{proof}

As a consequence, we derive the following corollary, which follows directly
from the last Theorem and Proposition \ref{CopoMethod}.

\begin{corollary}
The class $\mathcal{D}_{p}^{LipL}$ is the Banach Lip-Linear ideal generated
by the composition method from the Banach linear ideal $\mathcal{D}_{p}$. In
other words, 
\begin{equation*}
\mathcal{D}_{p}^{LipL}(X\times E;F)=\mathcal{D}_{p}\circ LipL_{0}(X\times
E;F),
\end{equation*}%
for all pointed metric space $X$ and Banach spaces $E,F.$
\end{corollary}

\section{Factorable strongly $p$-summing Lip-linear	operators}

The notion of factorability was introduced by Pellegrino et al. in \cite%
{FStr} in the context of multilinear operators, with the aim of identifying
the class of factorable strongly $p$-summing multilinear operators with the
class obtained by composition with absolutely $p$-summing linear operators.
Inspired by this idea, we extend the concept to the Lip-linear operators. We
begin by introducing the notion of strongly $p$-summing Lip-linear
operators, as a natural generalization of the concept introduced by Dimant
in \cite{Dim} for multilinear operators. We then define the corresponding
class of factorable operators within this framework. As an important
application, we establish a characterization of Hilbert spaces, originally
due to Kwapie\'{n} \cite{Kwa} in the linear case. Specifically, we show that
a Banach space $F$ is isomorphic to a Hilbert space if and only if every
factorable strongly $p$-summing Lip-linear operator with values in $F$ is
Cohen strongly $p$-summing.\textit{\ }Let\textit{\ }$p\in \left[ 1,\infty %
\right[ .$ Let $X$ be a pointed metric space and let $E,F$ be Banach spaces.
A Lip-linear operator $T:X\times E\longrightarrow F$\ is \textit{strongly }$p
$\textit{-summing}\ (in the sense of Dimant), if there exists a constant $C>0
$\ such that, for any $\left( x_{i}\right) _{i=1}^{n}\subset X$ and $\left(
e_{i}\right) _{i=1}^{n}\subset E$, we have%
\begin{equation}
(\sum\limits_{i=1}^{n}\left\Vert T(x_{i},e_{i})-T(y_{i},e_{i})\right\Vert
^{p})^{\frac{1}{p}}\leq C\sup_{\varphi \in B_{LipL_{0}\left( X\times
E\right) }}(\sum\limits_{i=1}^{n}\left\vert \varphi (x_{i},e_{i})-\varphi
(y_{i},e_{i})\right\vert ^{p})^{\frac{1}{p}}.  \label{4.1}
\end{equation}%
The class of all strongly $p$-summing Lip-linear operators from $X\times E$
into $F$, denoted by $\Pi _{St,p}^{LipL}(X\times E;F),$ forms a Banach space
equipped with the norm 
\begin{equation*}
\pi _{St,p}^{LipL}(T)=\inf \left\{ C:\text{ }C\text{ satisfies }\left( \ref%
{4.1}\right) \right\} .
\end{equation*}

Let us give an example. Let $X$ be a finite-dimensional Banach space. As
shown in \cite[Remark 3.3]{saadi}, the canonical Lipschitz embedding $\delta
_{X}:X\rightarrow \mathcal{F}\left( X\right) $ is Lipschitz $p$-summing with 
$\pi _{p}^{L}\left( \delta _{X}\right) =1$. Let $u:E\rightarrow F$ be a $p$%
-summing linear operator between Banach spaces. Define the Lip-linear
operator $T:X\times E\rightarrow X\widehat{\mathbb{\boxtimes }}_{\pi }F$ by 
\begin{equation*}
T\left( x,e\right) =\left( \delta _{X}\boxtimes u\right) \left( x,e\right)
=\delta _{x}\boxtimes u\left( e\right) .
\end{equation*}%
We claim that $T$ is strongly $p$-summing. Indeed, for any sequences $\left(
x_{i}\right) _{i=1}^{n},\left( y_{i}\right) _{i=1}^{n}\subset X$ and $\left(
e_{i}\right) _{i=1}^{n}\subset E$, we have%
\begin{eqnarray*}
(\sum\limits_{i=1}^{n}\left\Vert T(x_{i},e_{i})-T(y_{i},e_{i})\right\Vert
^{p})^{\frac{1}{p}} &=&(\sum\limits_{i=1}^{n}\left\Vert \left( \delta
_{x_{i}}-\delta _{y_{i}}\right) \boxtimes u\left( e_{i}\right) \right\Vert
^{p})^{\frac{1}{p}} \\
&=&(\sum\limits_{i=1}^{n}\left\Vert \delta _{x_{i}}-\delta
_{y_{i}}\right\Vert ^{p}\left\Vert u\left( e_{i}\right) \right\Vert ^{p})^{%
\frac{1}{p}}
\end{eqnarray*}%
Since $\delta _{X}$ is Lipschitz $p$-summing, we get%
\begin{eqnarray*}
(\sum\limits_{i=1}^{n}\left\Vert T(x_{i},e_{i})-T(y_{i},e_{i})\right\Vert
^{p})^{\frac{1}{p}} &\leq &\sup_{f\in
B_{X^{\#}}}(\sum\limits_{i=1}^{n}\left\Vert u\left( e_{i}\right)
\right\Vert ^{p}\left\vert f(x_{i})-f(y_{i})\right\vert ^{p})^{\frac{1}{p}}
\\
&\leq &\sup_{f\in B_{X^{\#}}}(\sum\limits_{i=1}^{n}\left\Vert u\left(
\left\vert f(x_{i})-f(y_{i})\right\vert e_{i}\right) \right\Vert ^{p})^{%
\frac{1}{p}}.
\end{eqnarray*}%
Using the fact that $u$ is $p$-summing, we obtain%
\begin{equation*}
(\sum\limits_{i=1}^{n}\left\Vert T(x_{i},e_{i})-T(y_{i},e_{i})\right\Vert
^{p})^{\frac{1}{p}}\leq \pi _{p}\left( u\right) \sup_{f\in
B_{X^{\#}},e^{\ast }\in B_{E^{\ast }}}(\sum\limits_{i=1}^{n}\left\vert
\left( f(x_{i})-f(y_{i})e^{\ast }\left( e_{i}\right) \right) \right\vert
^{p})^{\frac{1}{p}}.
\end{equation*}%
Thus,%
\begin{equation*}
(\sum\limits_{i=1}^{n}\left\Vert T(x_{i},e_{i})-T(y_{i},e_{i})\right\Vert
^{p})^{\frac{1}{p}}\leq \pi _{p}\left( u\right) \sup_{\varphi \in
B_{LipL_{0}\left( X\times E\right) }}(\sum\limits_{i=1}^{n}\left\vert
\varphi (x_{i},e_{i})-\varphi (y_{i},e_{i})\right\vert ^{p})^{\frac{1}{p}}.
\end{equation*}%
Hence, $T$ is strongly $p$-summing and satisfies 
\begin{equation*}
\pi _{St,p}^{LipL}(T)\leq \pi _{p}\left( u\right) .
\end{equation*}

An integral characterization of strongly $p$-summing Lip-linear operators is
obtained via the Pietsch Domination Theorem. The proof relies on a
generalized version of Pietsch's Theorem established by Pellegrino et al. in 
\cite{ps} and \cite{pss}.

\begin{theorem}
Let $1\leq p<\infty $ and $T\in LipL_{0}(X\times E;F).$ The following
assertions are equivalent:

1) $T:X\times E\longrightarrow F$\ is \textit{strongly }$p$\textit{-summing.}

2) There exist a regular probability measure $\mu $ on $B_{LipL_{0}\left(
X\times E\right) }$ and a constant $C>0$ such that for every $x,y\in X$ and $%
e\in E,$ we have%
\begin{equation*}
\left\Vert T(x,e)-T(y,e)\right\Vert \leq C(\int\limits_{B_{LipL_{0}\left(
X\times E\right) }}\left\vert \varphi (x,e)-\varphi (y,e)\right\vert
^{p}d\mu )^{\frac{1}{p}}.
\end{equation*}
\end{theorem}

\begin{proof}
Choosing the functions%
\begin{equation*}
\left\{ 
\begin{array}{l}
R:B_{LipL_{0}\left( X\times E\right) }\times \left( X\times X\times E\right)
\times \mathbb{R}\rightarrow \mathbb{R}^{+}:R\left( \varphi ,\left(
x,y,e\right) ,\lambda _{1}\right) =\left\vert \lambda _{1}\right\vert
\left\vert \varphi (x,e)-\varphi (y,e)\right\vert  \\ 
S:LipL_{0}(X,E;F)\times \left( X\times X\times E\right) \times \mathbb{R}%
\times \mathbb{R}\rightarrow \mathbb{R}^{+}: \\ 
S\left( T,\left( x,y,e\right) ,\lambda _{1},\lambda _{2}\right) =\left\vert
\lambda _{1}\right\vert \left\Vert T\left( x,e\right) -T\left( y,e\right)
\right\Vert .%
\end{array}%
\right. 
\end{equation*}%
These maps satisfy conditions $\left( 1\right) $ and $\left( 2\right) $ from 
\cite[Page 1255]{pss}, allowing us to conclude that $T:X\times E\rightarrow F
$ is strongly $p$-summing if and only if 
\begin{equation*}
(\sum_{i=1}^{n}\left\Vert S\left( T,\left( x_{i},y_{i},e_{i}\right) ,\lambda
_{i,1},\lambda _{i,2}\right) \right\Vert ^{p})^{\frac{1}{p}}\leq
\sup_{\varphi \in B_{LipL_{0}\left( X\times E\right)
}}(\sum_{i=1}^{n}R\left( \varphi ,\left( x_{i},y_{i},e_{i}\right) ,\lambda
_{i,1}\right) ^{p})^{\frac{1}{p}}.
\end{equation*}%
Thus, $T$ is $R$-$S$-abstract $p$-summing. A result in \cite[Theorem 4.6]%
{pss} states that $T$ is $R$-$S$-abstract $p$-summing if and only if there
exists a positive constant $C$ and Radon probability measures $\mu $ on $%
B_{LipL_{0}\left( X\times E\right) }$such that 
\begin{equation*}
S\left( T,\left( x,y,e\right) ,\lambda _{1},\lambda _{2}\right) \leq
C(\int_{B_{LipL\left( X\times E\right) }}R\left( \varphi ,\left(
x,y,e\right) ,\lambda _{1}\right) ^{p}d\mu \left( \varphi \right) )^{\frac{1%
}{p}}.
\end{equation*}%
Consequently, 
\begin{equation*}
\left\Vert T(x,e)-T(y,e)\right\Vert \leq C(\int\limits_{B_{LipL_{0}\left(
X\times E\right) }}\left\vert \varphi (x,e)-\varphi (y,e)\right\vert
^{p}d\mu \left( \varphi \right) )^{\frac{1}{p}},
\end{equation*}%
this concludes the proof.
\end{proof}

\begin{proposition}
Let $T\in LipL_{0}(X\times E;F)$ such that its linearization $\widehat{T}:X%
\widehat{\boxtimes }_{\pi }E\longrightarrow F$\ is $p$\textit{-summing, then 
}$T$ is strongly $p$-summing\textit{.}
\end{proposition}

\begin{proof}
Suppose that $\widehat{T}$ is $p$-summing. Let $\left( x_{i}\right)
_{i=1}^{n}\subset X$ and $\left( e_{i}\right) _{i=1}^{n}\subset E$, we have%
\begin{eqnarray*}
\sum\limits_{i=1}^{n}\left\Vert T(x_{i},e_{i})-T(y_{i},e_{i})\right\Vert
^{p})^{\frac{1}{p}} &=&\sum\limits_{i=1}^{n}\left\Vert \widehat{T}(\delta
_{\left( x_{i},y_{i}\right) }\boxtimes e_{i})\right\Vert ^{p})^{\frac{1}{p}}
\\
&\leq &\pi _{p}\left( \widehat{T}\right) \sup_{\widehat{\varphi }\in
B_{\left( X\widehat{\boxtimes }_{\pi }E\right) ^{\ast
}}}(\sum\limits_{i=1}^{n}\left\vert \widehat{\varphi }(\delta _{\left(
x_{i},y_{i}\right) }\boxtimes e_{i})\right\vert ^{p})^{\frac{1}{p}} \\
&\leq &\pi _{p}\left( \widehat{T}\right) \sup_{\varphi \in B_{LipL_{0}\left(
X\times E\right) }}(\sum\limits_{i=1}^{n}\left\vert \varphi
(x_{i},e_{i})-\varphi (y_{i},e_{i})\right\vert ^{p})^{\frac{1}{p}}.
\end{eqnarray*}%
Hence, $T$ is strongly $p$-summing and 
\begin{equation*}
\pi _{St,p}^{LipL}(T)\leq \pi _{p}\left( \widehat{T}\right) .
\end{equation*}
\end{proof}

The converse of the previous Proposition does not hold. Consider the
following counterexample. Let $\sigma ^{LipL}:\mathbb{R}\times \mathbb{%
R\rightarrow R}\widehat{\mathbb{\boxtimes }}_{\pi }\mathbb{R}\mathbb{\ }$be
defined by%
\begin{equation*}
\sigma ^{LipL}=\delta _{\mathbb{R}}\boxtimes id_{\mathbb{R}}.
\end{equation*}%
Since $id_{\mathbb{R}}$ is $p$-summing (as $\mathbb{R}$ is
finite-dimensional), it follows that $\sigma ^{LipL}$ is strongly $p$%
-summing. However, observe the following commutative diagram%
\begin{equation*}
\begin{array}{ccc}
\mathbb{R}\times \mathbb{R} & \overset{\sigma ^{LipL}}{\longrightarrow } & 
\mathbb{R}\widehat{\mathbb{\boxtimes }}_{\pi }\mathbb{R} \\ 
\sigma ^{LipL}\downarrow  & \nearrow \widehat{\sigma ^{LipL}} &  \\ 
\mathbb{R}\widehat{\mathbb{\boxtimes }}_{\pi }\mathbb{R} &  & 
\end{array}%
\end{equation*}%
In this setting, we have $\widehat{\sigma ^{LipL}}=id_{\mathbb{R}\widehat{%
\mathbb{\boxtimes }}_{\pi }\mathbb{R}}$, which is not $p$-summing. This
follows from the fact that $\mathbb{R}\widehat{\mathbb{\boxtimes }}_{\pi }%
\mathbb{R}$ is infinite-dimensional (in fact, it is isometrically isomorphic
to $\mathcal{F}\left( \mathbb{R}\right) $)$.$ This observation motivates the
introduction of the class of factorable strongly $p$-summing Lip-linear
operators. This new class generalizes the notion of factorable strongly $p$%
-summing multilinear operators introduced by Pellegrino et al. in \cite{FStr}%
.

\begin{definition}
A Lip-linear operator $T:X\times E\longrightarrow F$\ is \textit{factorable} 
\textit{strongly }$p$\textit{-summing}\ if there exists a constant $C>0$\
such that, for any $\left( x_{i}^{j}\right) _{i=1}^{n},\left(
y_{i}^{j}\right) _{i=1}^{n}\subset X,\left( 1\leq j\leq m\right) ,\left(
e_{i}\right) _{i=1}^{n}\subset E$, we have%
\begin{equation}
\sum\limits_{i=1}^{n}(\left\Vert
\sum\limits_{j=1}^{m}T(x_{i}^{j},e_{i}^{j})-T(y_{i}^{j},e_{i}^{j})\right%
\Vert ^{p})^{\frac{1}{p}}\leq C\sup_{\varphi \in B_{LipL_{0}\left( X\times
E\right) }}(\sum\limits_{i=1}^{n}\left\vert \sum\limits_{j=1}^{m}\varphi
(x_{i}^{j},e_{i}^{j})-\varphi (y_{i}^{j},e_{i}^{j})\right\vert ^{p})^{\frac{1%
}{p}}.  \label{4.2}
\end{equation}%
The class of all factorable strongly $p$-summing Lip-linear operators from $%
X\times E$ into $F$, which is denoted by $\Pi _{F,St,p}^{LipL}\left( X\times
E;F\right) ,$ is a Banach space with the norm 
\begin{equation*}
\pi _{F,St,p}^{LipL}(T)=\inf \left\{ C:\text{ }C\text{ satisfies }\left( \ref%
{4.2}\right) \right\} .
\end{equation*}
\end{definition}

By setting $m=1$ in the formula $\left( \ref{4.2}\right) $, we recover the
definition of strongly $p$-summing Lipschitz-linear operators. That is, 
\begin{equation*}
\Pi _{F,St,p}^{LipL}\left( X\times E;F\right) \subset \Pi
_{St,p}^{LipL}(X\times E;F).
\end{equation*}

We now establish the connection between a Lip-linear operator and its
linearization in the context of $p$-summing operators. This correspondence
highlights an important aspect that is absent from the definition of
strongly $p$-summing operators. Moreover, this relation allows us to derive
the Pietsch Domination Theorem for the class of factorable strongly $p$%
-summing Lip-linear operators.

\begin{theorem}
\label{FactLin}\textit{Let} $1\leq p<\infty .$ Let $T:X\times E\rightarrow F$
be a Lip-linear operator. \textit{The following properties are equivalent.}

1) \textit{The operator }$T$ \textit{belongs to }$\Pi _{F,St,p}^{LipL}\left(
X\times E;F\right) .$

2) \textit{The linearization operator }$\widehat{T}$\textit{\ is }$p$\textit{%
-summing.}

3) There exist a regular probability measure $\mu $ on $B_{LipL_{0}\left(
X\times E\right) }$ and a constant $C>0$ such that for every $\left(
x^{j}\right) _{j=1}^{m},\left( y^{j}\right) _{j=1}^{m}\subset X$ and $\left(
e^{j}\right) _{j=1}^{m}\subset E,$ we have%
\begin{equation}
\left\Vert \sum\limits_{j=1}^{m}T(x^{j},e^{j})-T(y^{j},e^{j})\right\Vert
\leq C(\int\limits_{B_{LipL_{0}\left( X\times E\right) }}\left\vert
\sum\limits_{j=1}^{m}\varphi (x^{j},e^{j})-\varphi (y^{j},e^{j})\right\vert
^{p}d\mu )^{\frac{1}{p}}  \label{4.3}
\end{equation}%
\textit{In add, we have the following isometrically identification}%
\begin{equation*}
\Pi _{F,St,p}^{LipL}\left( X\times E;F\right) =\Pi _{p}\left( X\widehat{%
\mathbb{\boxtimes }}_{\pi }E;F\right) \text{\textit{.}}
\end{equation*}
\end{theorem}

\begin{proof}
$1)\Longrightarrow 2):$ Let $T\in \Pi _{F,St,p}^{LipL}\left( X\times
E;F\right) $ and $u_{i}=\sum\limits_{j=1}^{m}\delta _{\left(
x_{i}^{j},y_{i}^{j}\right) }\boxtimes e_{i}^{j}\in X\widehat{\mathbb{%
\boxtimes }}_{\pi }E\left( 1\leq i\leq n\right) .$ Then%
\begin{eqnarray*}
\left( \sum\limits_{i=1}^{n}\left\Vert \widehat{T}\left( u_{i}\right)
\right\Vert ^{p}\right) ^{\frac{1}{p}} &=&(\sum\limits_{i=1}^{n}\left\Vert
\sum\limits_{j=1}^{m}T(x_{i}^{j},e_{i}^{j})-T(y_{i}^{j},e_{i}^{j})\right%
\Vert ^{p})^{\frac{1}{p}} \\
&\leq &\pi _{F,St,p}^{LipL}\left( T\right) \sup_{\varphi \in
B_{LipL_{0}\left( X\times E\right) }}(\sum\limits_{i=1}^{n}\left\vert
\sum\limits_{j=1}^{m}\varphi (x_{i}^{j},e_{i}^{j})-\varphi
(y_{i}^{j},e_{i}^{j})\right\vert ^{p})^{\frac{1}{p}} \\
&\leq &\pi _{F,St,p}^{LipL}\left( T\right) \sup_{\varphi \in
B_{LipL_{0}\left( X\times E\right) }}(\sum\limits_{i=1}^{n}\left\vert 
\widehat{\varphi }\left( u_{i}\right) \right\vert ^{p})^{\frac{1}{p}} \\
&\leq &\pi _{F,St,p}^{LipL}\left( T\right) \sup_{\widehat{\varphi }\in
\left( X\widehat{\mathbb{\boxtimes }}_{\pi }E\right) ^{\ast
}}(\sum\limits_{i=1}^{n}\left\vert \widehat{\varphi }\left( u_{i}\right)
\right\vert ^{p})^{\frac{1}{p}}.
\end{eqnarray*}%
Hence$,$ $\widehat{T}$ is $p$-summing and 
\begin{equation*}
\pi \left( \widehat{T}\right) \leq \pi _{F,St,p}^{LipL}\left( T\right) .
\end{equation*}%
$2)\Longrightarrow 3):$ Suppose that $\widehat{T}$ is $p$-summing. We now
proceed to prove $\left( \ref{4.3}\right) $. By \cite[Theorem 2.12]{distel},
there exists a Radon probability measures $\mu $ on $B_{\left( X\widehat{%
\mathbb{\boxtimes }}_{\pi }E\right) ^{\ast }}$ such that%
\begin{equation*}
\left\Vert \widehat{T}\left( u\right) \right\Vert \leq \pi \left( \widehat{T}%
\right) (\int\limits_{B_{\left( X\widehat{\mathbb{\boxtimes }}_{\pi
}E\right) ^{\ast }}}\left\vert \varphi (u)\right\vert ^{p}d\mu \left(
\varphi \right) )^{\frac{1}{p}}.
\end{equation*}%
Let $\left( x^{j}\right) _{j=1}^{m},\left( y^{j}\right) _{j=1}^{m}\subset X$
and $\left( e^{j}\right) _{j=1}^{m}\subset E.$ We put $u=\sum%
\limits_{j=1}^{m}\delta _{\left( x^{j},y^{j}\right) }\boxtimes e^{j}$, we
btain%
\begin{eqnarray*}
\left\Vert \sum\limits_{j=1}^{m}T(x^{j},e^{j})-T(y^{j},e^{j})\right\Vert 
&=&\left\Vert \widehat{T}\left( u\right) \right\Vert  \\
&\leq &\pi _{p}\left( \widehat{T}\right) (\int\limits_{B_{\left( X\widehat{%
\mathbb{\boxtimes }}_{\pi }E\right) ^{\ast }}}\left\vert \varphi
(u)\right\vert ^{p}d\mu )^{\frac{1}{p}} \\
&\leq &\pi _{p}\left( \widehat{T}\right) (\int\limits_{B_{LipL_{0}\left(
X\times E\right) }}\left\vert \sum\limits_{j=1}^{m}\varphi
(x^{j},e^{j})-\varphi (y^{j},e^{j})\right\vert ^{p}d\mu )^{\frac{1}{p}}.
\end{eqnarray*}%
$3)\Longrightarrow 1):$ Let $\left( x_{i}^{j}\right) _{i=1}^{n},\left(
y_{i}^{j}\right) _{i=1}^{n}\subset X$ and $\left( e_{i}^{j}\right)
_{i=1}^{n}\subset E\left( 1\leq j\leq m\right) $. For $1\leq i\leq n,$ we
have

$\left\Vert
\sum\limits_{j=1}^{m}T(x_{i}^{j},e_{i}^{j})-T(y_{i}^{j},e_{i}^{j})\right%
\Vert ^{p}\leq C^{p}\int\limits_{B_{LipL_{0}\left( X\times E\right)
}}\left\vert \sum\limits_{j=1}^{m}\varphi (x_{i}^{j},e_{i}^{j})-\varphi
(y_{i}^{j},e_{i}^{j})\right\vert ^{p}d\mu .$

Then

$\sum\limits_{i=1}^{n}\left\Vert
\sum\limits_{j=1}^{m}T(x_{i}^{j},e_{i}^{j})-T(y_{i}^{j},e_{i}^{j})\right%
\Vert ^{p}$

$\leq C^{p}\sum\limits_{i=1}^{n}\int\limits_{B_{LipL_{0}\left( X\times
E\right) }}\left\vert \sum\limits_{j=1}^{m}\varphi
(x_{i}^{j},e_{i}^{j})-\varphi (y_{i}^{j},e_{i}^{j})\right\vert ^{p}d\mu .$

$\leq C^{p}\int\limits_{B_{LipL_{0}\left( X\times E\right)
}}\sum\limits_{i=1}^{n}\left\vert \sum\limits_{j=1}^{m}\varphi
(x_{i}^{j},e_{i}^{j})-\varphi (y_{i}^{j},e_{i}^{j})\right\vert ^{p}d\mu .$

Taking the supremum over all $\varphi \in B_{LipL_{0}\left( X\times E\right)
},$ we obtain

$\sum\limits_{i=1}^{n}\left\Vert
\sum\limits_{j=1}^{m}T(x_{i}^{j},e_{i}^{j})-T(y_{i}^{j},e_{i}^{j})\right%
\Vert ^{p}\leq C^{p}\sup_{\varphi \in B_{LipL_{0}\left( X\times E\right)
}}\sum\limits_{i=1}^{n}\left\vert \sum\limits_{j=1}^{m}\varphi
(x_{i}^{j},e_{i}^{j})-\varphi (y_{i}^{j},e_{i}^{j})\right\vert ^{p}$.

Hence, $T$ is factorable strngly $p$-summing.
\end{proof}

As a consequence, we derive the following corollary, which follows directly
from the above Theorem and Proposition $\ref{CopoMethod}$.

\begin{corollary}
The class $\Pi _{F,St,p}^{LipL}$ is the Banach Lip-Linear ideal generated by
the composition method from the Banach linear ideal $\Pi _{p}$. In other
words, 
\begin{equation*}
\Pi _{F,St,p}^{LipL}(X\times E;F)=\Pi _{p}\circ LipL_{0}(X\times E;F),
\end{equation*}%
for all pointed metric space $X$ and Banach spaces $E,F.$
\end{corollary}

We now present an important characterization of Hilbert spaces using the
concepts of factorable strongly $2$-summing and Cohen strongly $2$-summing
Lip-linear operators.

\begin{theorem}
Let $F$ be a Banach space. The following assertions are equivalent.

1) $F$ is isomorphic to a Hilbert space.

2) For every pointed metric space $X$ and Banach space $E$, we have%
\begin{equation*}
\Pi _{F,St,2}^{LipL}(X\times E;F)\subset \mathcal{D}_{2}^{LipL}(X\times E;F).
\end{equation*}
\end{theorem}

\begin{proof}
$1)\Longrightarrow 2):$ Assume that $F$ is isomorphic to a Hilbert space. By
a classical result of Kwapie\'{n} \cite{Kwa}, for every Banach space $G$,
every $2$-summing linear operator $u:G\rightarrow F$ is also Cohen strongly $%
2$-summing, i.e.,%
\begin{equation*}
\Pi _{2}\left( G;F\right) \subset \mathcal{D}_{2}\left( G;F\right) .
\end{equation*}%
Now, let $T\in \Pi _{F,St,2}^{LipL}(X\times E;F)$ where $X$ is a pointed
metric space and $E$ is a Banach space. By Theorem $\ref{FactLin}$, the
linearization $\widehat{T}:X\widehat{\mathbb{\boxtimes }}_{\pi }E\rightarrow
F$ is $2$-summing. Thus, by Kwapie\'{n}'s result, $\widehat{T}$ is strongly $%
2$-summing. Consequently, by Theorem $\ref{T_L}$, $T\in \mathcal{D}%
_{2}^{LipL}(X\times E;F)$.

$2)\Longrightarrow 1):$ Let $E$ be any Banach space and $u:E\rightarrow F$ a 
$2$-summing linear operator. We aim to prove that $u$ is strongly $2$%
-summing. To do this, define the Lip-linear (bilinear) operator $T:\mathbb{%
R\times }E\rightarrow F$ by%
\begin{equation*}
T\left( \lambda ,e\right) =\lambda u\left( e\right) .
\end{equation*}%
We show that $T\in \Pi _{F,St,2}^{LipL}(\mathbb{R}\times E;F).$ Let $\left(
\lambda _{i}^{j}\right) _{i=1}^{n},\left( \gamma _{i}^{j}\right)
_{i=1}^{n}\subset \mathbb{R}$ and $\left( e_{i}^{j}\right) _{i=1}^{n}\subset
E\left( 1\leq j\leq m\right) .$ Then%
\begin{eqnarray*}
\sum\limits_{i=1}^{n}(\left\Vert \sum\limits_{j=1}^{m}T(\lambda
_{i}^{j},e_{i}^{j})-T(\gamma _{i}^{j},e_{i}^{j})\right\Vert ^{2})^{\frac{1}{2%
}} &=&\sum\limits_{i=1}^{n}(\left\Vert \sum\limits_{j=1}^{m}\left( \lambda
_{i}^{j}-\gamma _{i}^{j}\right) u\left( e_{i}^{j}\right) \right\Vert ^{2})^{%
\frac{1}{2}} \\
&=&\sum\limits_{i=1}^{n}(\left\Vert u\left( \sum\limits_{j=1}^{m}\left(
\lambda _{i}^{j}-\gamma _{i}^{j}\right) e_{i}^{j}\right) \right\Vert ^{2})^{%
\frac{1}{2}}
\end{eqnarray*}%
Since $u$ is $2$-summing,\ we have%
\begin{eqnarray*}
&&\sum\limits_{i=1}^{n}(\left\Vert \sum\limits_{j=1}^{m}T(\lambda
_{i}^{j},e_{i}^{j})-T(\gamma _{i}^{j},e_{i}^{j})\right\Vert ^{2})^{\frac{1}{2%
}} \\
&\leq &\pi _{2}\left( u\right) \sup_{e^{\ast }\in B_{E^{\ast
}}}\sum\limits_{i=1}^{n}(\left\vert \sum\limits_{j=1}^{m}\lambda
_{i}^{j}e^{\ast }\left( e_{i}^{j}\right) -\gamma _{i}^{j}e^{\ast }\left(
e_{i}^{j}\right) \right\vert ^{2})^{\frac{1}{2}} \\
&\leq &\pi _{2}\left( u\right) \sup_{e^{\ast }\in B_{E^{\ast
}}}\sum\limits_{i=1}^{n}(\left\vert \sum\limits_{j=1}^{m}id_{\mathbb{R}%
}\left( \lambda _{i}^{j}\right) e^{\ast }\left( e_{i}^{j}\right) -id_{%
\mathbb{R}}\left( \gamma _{i}^{j}\right) e^{\ast }\left( e_{i}^{j}\right)
\right\vert ^{2})^{\frac{1}{2}} \\
&\leq &\pi _{2}\left( u\right) \sup_{\varphi \in B_{LipL_{0}\left( X\times
E\right) }}(\sum\limits_{i=1}^{n}\left\vert \sum\limits_{j=1}^{m}\varphi
(\lambda _{i}^{j},e_{i}^{j})-\varphi (\gamma _{i}^{j},e_{i}^{j})\right\vert
^{2})^{\frac{1}{2}}.
\end{eqnarray*}%
Thus, $T\in \Pi _{F,St,2}^{LipL}(\mathbb{R}\times E;F).$ By assumption $2)$, 
$T\in \mathcal{D}_{2}^{LipL}(\mathbb{R}\times E;F)$. Now, let $\left(
e_{i}\right) _{i=1}^{n}\subset E$ and $\left( z_{i}^{\ast }\right)
_{i=1}^{n}\subset F^{\ast }.$ Define $\lambda _{i}=1$ and $\gamma _{i}=0$ $%
\left( 1\leq i\leq n\right) .$ Then, 
\begin{eqnarray*}
\sum\limits_{i=1}^{n}\left\vert \left\langle u\left( e_{i}\right)
,z_{i}^{\ast }\right\rangle \right\vert  &=&\sum\limits_{i=1}^{n}\left\vert
\left\langle \lambda _{i}u\left( e_{i}\right) -\gamma _{i}\times u\left(
e_{i}\right) ,z_{i}^{\ast }\right\rangle \right\vert  \\
&=&\sum\limits_{i=1}^{n}\left\vert \left\langle T\left( 1,e_{i}\right)
-T\left( 0,e_{i}\right) ,z_{i}^{\ast }\right\rangle \right\vert ,
\end{eqnarray*}%
Using the definition of Cohen strongly $2$-summing Lip-linear operators, we
get%
\begin{eqnarray*}
\sum\limits_{i=1}^{n}\left\vert \left\langle u\left( e_{i}\right)
,z_{i}^{\ast }\right\rangle \right\vert  &\leq &\pi _{F,St,2}^{LipL}\left(
T\right) (\sum\limits_{i=1}^{n}\left( d(1,0)\left\Vert e_{i}\right\Vert
\right) ^{2})^{\frac{1}{2}}\sup_{z^{\ast \ast }\in B_{F^{\ast \ast
}}}(\sum\limits_{i=1}^{n}\left\vert z^{\ast \ast }(z_{i}^{\ast
})\right\vert ^{2})^{\frac{1}{2}} \\
&\leq &\pi _{F,St,2}^{LipL}\left( T\right) (\sum\limits_{i=1}^{n}\left\Vert
e_{i}\right\Vert ^{2})^{\frac{1}{2}}\sup_{z^{\ast \ast }\in B_{F^{\ast \ast
}}}(\sum\limits_{i=1}^{n}\left\vert z^{\ast \ast }(z_{i}^{\ast
})\right\vert ^{2})^{\frac{1}{2}}.
\end{eqnarray*}%
Therefore, $u$ is strongly $2$-summing. By Kwapie\'{n}'s theorem, this
implies that $F$ is isomorphic to a Hilbert space.
\end{proof}

\end{document}